\pdfoutput=1
\RequirePackage{ifpdf}
\ifpdf 
\documentclass[pdftex]{sigma}
\else
\documentclass{sigma}
\fi

\usepackage{mathrsfs}
\usepackage[all]{xy}

\numberwithin{equation}{section}
\newtheorem{Theorem}{Theorem}[section]
\newtheorem{Corollary}[Theorem]{Corollary}
\newtheorem{Lemma}[Theorem]{Lemma}
{ \theoremstyle{definition}
\newtheorem{Remark}[Theorem]{Remark} }

\DeclareMathOperator{\ad}{ad}
\DeclareMathOperator{\co}{co}
\DeclareMathOperator{\Tor}{Tor}

\newcommand{\Exterior}{\mathchoice{{\textstyle\bigwedge}}%
    {{\bigwedge}}%
    {{\textstyle\wedge}}%
    {{\scriptstyle\wedge}}}

\newcommand{\ps}[1]{~\hspace{-4pt}_{^{(#1)}}}

\newcommand{\ns}[1]{~\hspace{-4pt}_{_{{\langle#1\rangle}}}}

\newcommand{\wbar}[1]{\overline{#1}}

\begin{document}

\newcommand{\arXivNumber}{1412.2359}

\allowdisplaybreaks

\renewcommand{\PaperNumber}{041}

\FirstPageHeading

\ShortArticleName{Cyclic Homology and Quantum Orbits}

\ArticleName{Cyclic Homology and Quantum Orbits}

\Author{Tomasz MASZCZYK and Serkan S\"UTL\"U}

\AuthorNameForHeading{T.~Maszczyk and S.~S\"utl\"u}

\Address{Institute of Mathematics, University of Warsaw, Warsaw, Poland}
\Email{\href{mailto:t.maszczyk@uw.edu.pl}{t.maszczyk@uw.edu.pl},
\href{mailto:serkansutlu@gmail.com}{serkansutlu@gmail.com}}

\ArticleDates{Received December 07, 2014, in f\/inal form May 12, 2015; Published online May 30, 2015}

\Abstract{A natural isomorphism between the cyclic object computing the relative cyclic homology of a~homogeneous
quotient-coalgebra-Galois extension, and the cyclic object computing the cyclic homology of a~Galois coalgebra with SAYD
coef\/f\/icients is presented.
The isomorphism can be viewed as the cyclic-homological counterpart of the Takeuchi--Galois correspondence between the
left coideal subalgebras and the quotient right module coalgebras of a~Hopf algebra.
A~spectral sequence generalizing the classical computation of Hochschild homology of a~Hopf algebra to the case of
arbitrary homogeneous quotient-coalgebra-Galois extensions is constructed.
A~Pontryagin type self-duality of the Takeuchi--Galois correspondence is combined with the cyclic duality of Connes in
order to obtain dual results on the invariant cyclic homology, with SAYD coef\/f\/icients, of algebras of invariants in
homogeneous quotient-coalgebra-Galois extensions.
The relation of this dual result with the Chern cha\-rac\-ter, Frobenius reciprocity, and inertia phenomena in the local
Langlands program, the Chen--Ruan--Brylinski--Nistor orbifold cohomology and the Clif\/ford theory is discussed.}

\Keywords{cyclic homology; homogenous quotient-coalgebra-Galois extensions; Takeuchi--Galois correspondence; Pontryagin
duality}

\Classification{19D55; 57T15; 06A15; 46A20}

\section{Introduction}

It is well known due to Takeuchi~\cite{Take72} (vastly extended by van Oystaeyen--Zhang~\cite{VanOZhan94} and
Schauenburg~\cite{Scha98-II}) that given a~Hopf algebra ${\cal H}$ there is a~Galois 1-1 correspondence between its left
coideal subalgebras, for which ${\cal H}$ is a~faithfully f\/lat algebra extension, and its quotient right module
coalgebras, for which ${\cal H}$ is a~faithfully cof\/lat coalgebra coextension.

Such a~correspondence can also be viewed as a~relation between the extensions of comodule algebras and the coextensions
of module coalgebras, both of which having their specif\/ic homological invariants (Hochschild, cyclic, periodic cyclic
and negative cyclic homology) computed from appropriate cyclic objects.

For algebra extensions, it is a~relative cyclic object introduced by Kadison~\cite{Kadi89}, and for comodule algebras it
is a~cyclic object with stable anti-Yetter--Drinfeld (SAYD) coef\/f\/icients introduced~by
Hajac--Khalkhali--Rangipour--Sommerh\"{a}user~\cite{HajaKhalRangSomm04-II,HajaKhalRangSomm04-I}, and independently by
Jara--\c{S}tefan~\cite{JaraStef06} (in a~cyclic dual version).

For module coalgebras, it is a~cyclic object with SAYD coef\/f\/icients which is cyclic dual to that
of~\cite{HajaKhalRangSomm04-II}, and for coalgebra extensions it is the cyclic dual of the Pontryagin dual analogue
of~\cite{Kadi89}.

Therefore it is very natural to ask whether these dif\/ferent types of cyclic objects are related in the context of the
aforementioned Takeuchi--Galois correspondence.

This question goes far beyond the Galois theory (herein the theory of so called homogeneous quotient coalgebra-Galois
extensions~\cite{BrzeHaja09}) and reaches topology.
In the case of smooth functions on compact Lie groups and their homogeneous spaces, relative periodic cyclic homology
computes the vector bundle of de Rham cohomology of the stabilizer over the homogeneous space.
This bundle is equipped with the Gauss--Manin connection (see~\cite{Getz91-92} for noncommutative f\/ibrations with
commutative base) determining a~local system of coef\/f\/icients whose cohomology appears in the second page of the Leray
spectral sequence interpolating between the cohomology of the compact Lie group and the cohomology of the homogeneous space.
Then the above Takeuchi--Galois correspondence boils down to one between stabilizers and orbits, and the algebra
extension describes the orbital map.

Moreover, the cyclic object, with SAYD coef\/f\/icients, of a~module coalgebra is cyclic dual to the cocyclic object
generalizing the one used by Connes--Moscovici in the proof of a~gene\-ra\-li\-zed transversal local index theorem for
foliated spaces, where the index computation relies on a~symmetry governed by a~Hopf
algebra~\cite{ConnMosc98,ConnMosc99,HajaKhalRangSomm04-II}.

The aim of the present paper is to prove that these two types of homological invariants are isomorphic in a~natural way.
More explicitly,

\smallskip

\noindent\textbf{Theorem~3.4.} \emph{Let ${\cal I}\subseteq {\cal H}$ be a~coideal right ideal in a~Hopf algebra ${\cal H}$
such that ${\cal H}^{\co{\cal H}/{\cal I}}\subseteq {\cal H}$ is an ${\cal H}/{\cal I}$-Galois extension.
Let also $\ad({\cal H})={\cal H}$ be the left-right SAYD module with the right adjoint action, and the left
coaction given by the comultiplication of ${\cal H}$.
Then there exists an isomorphism
\begin{gather*}
\psi_{n}\colon \
 {\rm C}_{n}({\cal H}/{\cal I},\ad({\cal H}))_{{\cal H}}\longrightarrow {\rm C}_{n}\big({\cal H}\mid{\cal H}^{\co{\cal H}/{\cal I}}\big)
\end{gather*}
of cyclic modules, defined by~\eqref{aux-map-psi}, \eqref{aux-map-vp}.}

\smallskip

Since the domain of the isomorphism depends only on the right ${\cal H}$-module quotient coalgebra of ${\cal H}$, and
the codomain depends on the left comodule subalgebra of ${\cal H}$, this isomorphism can be viewed as
a~\emph{cyclic-homological Takeuchi--Galois transform} accompanying the Takeuchi--Galois
transform ${\cal H}/{\cal I}\mapsto {\cal H}^{\co{\cal H}/{\cal I}}$.

This perspective manifested itself in the paper of Jara--\c{S}tefan~\cite{JaraStef06}, where the f\/irst example of the
cyclic homology with SAYD coef\/f\/icients dif\/ferent from the modular pairs in involution of Connes--Moscovici was presented.
Their main result~\cite[Theorem~3.7]{JaraStef06} is an isomorphism of cyclic objects computing relative cyclic homology
of a~Hopf--Galois extension, and cyclic homology of the Galois Hopf algebra itself with an appropriate SAYD coef\/f\/icients
(def\/ined via the Miyashita--Ulbrich action of the Hopf algebra).

However, the Hopf--Galois context is too narrow to cover the full context of the Takeuchi--Galois correspondence (making
sense only for a~restricted correspondence between left comodule subalgebras ${\cal B}\subseteq {\cal H}$ such that ${\cal B}^{+}{\cal H}$ is a~Hopf ideal and quotient Hopf-algebras~\cite{Take72}) and misses
important examples such as non-standard Podle\'s quantum spheres.
In such cases, instead of a~Hopf algebra of Galois symmetry one has merely a~coaugmented quotient coalgebra of a~bigger
Hopf-algebra acting transitively on both the base and the total quantum space of a~quantum principal bundle (so called
\emph{homogeneous quotient-coalgebra-Galois extension}).
Then a~Miyashita--Ulbrich type action doesn't make sense, and the Jara--\c{S}tefan isomorphism~\cite[Theorem~3.7]{JaraStef06} cannot be applied.

Nevertheless, as we show below, when restricted to the homogeneous Hopf-quotient-Galois case our isomorphism and that of
Jara--\c{S}tefan~\cite[Theorem~3.7]{JaraStef06}, despite the apparently dif\/ferent def\/initions, essentially coincide.

A deeper conceptual motivation of the coalgebra-Galois approach, beyond the context of the Takeuchi--Galois
correspondence, comes from the so called \emph{coisotropic creed} of Poisson geo\-met\-ry~\cite{Lu93}.
The latter, accepted to bypass the problem of poverty of Poisson subgroups of Poisson groups, launched the theory of
coisotropic subgroups.
The quantized counterpart of coisotropic subgroups of a~Poisson group is a~quotient coaugmented right module coalgebra
of a~Hopf algebra by a~coideal right ideal.
The non-standard Podle\'s quantum spheres are examples of quantum orbits corresponding to such coalgebraic quantum
stabilizers~\cite{DijkKoor94}.

It is worth noticing that even in such a~simple homogeneous ${\cal H}$-Galois extension as $k={\cal H}^{\co{\cal H}}\to{\cal H}$,
computing the relative Hochschild homology (being the Hochschild homology of the Hopf algebra itself)
is nontrivial and of fundamental interest~\cite{Bich13,BrowZhan08,CollHartThom09,FengTsyg91,GinzKuma93,HadfKrah05}.
The following folklore result, proven explicitly f\/irst by Feng--Tsygan~\cite{FengTsyg91}, next by Bichon~\cite{Bich13} by
a~similar method with a~dif\/ferent perspective, but in fact already implicitly contained in Cartan--Eilenberg's book (it
is enough to dualize a~consequence of~\cite[Theorem~VIII.3.1]{CartEile56} and rephrase the result in terms of Hopf
algebras and Hopf algebra homology) says that

\smallskip

\noindent\textbf{Corollary~3.6.} \emph{For any Hopf algebra ${\cal H}$,
\begin{gather*}
{\rm HH}_\bullet({\cal H})=\Tor^{{\cal H}}_{\bullet}(k, \ad({\cal H})),
\end{gather*}
where on the right hand side the left ${\cal H}$-module structure on ${\cal H}$ comes from its canonical SAYD module structure.}

\smallskip

We generalize this result to the case of arbitrary \emph{homogeneous quotient coalgebra-Galois extensions} as follows

\smallskip

\noindent\textbf{Theorem~3.5.} \emph{Let ${\cal I}\subseteq {\cal H}$ be a~coideal right ideal in a~Hopf algebra ${\cal H}$
such that ${\cal H}^{\co{\cal H}/{\cal I}}\subseteq {\cal H}$ be a~homogeneous ${\cal H}/{\cal I}$-Galois extension.
Then there exists a~spectral sequence $($constructed in the proof$)$ such that
\begin{gather*}
{\rm HH}_\bullet\big({\cal H}\mid {\cal H}^{\co{\cal H}/{\cal I}}\big)={\rm E}^2_{\bullet,0},
\qquad
{\rm E}^2_{\bullet,\bullet}\Longrightarrow \Tor^{{\cal H}}_{\bullet}(k, \ad({\cal H})).
\end{gather*}
In particular, we have a~five-term exact sequence}
\begin{gather*}
\phantom{\dots\to }
\Tor^{{\cal H}}_2(k,\ad({\cal H}))\to {\rm HH}_2\big({\cal H}\mid {\cal H}^{\co{\cal H}/{\cal I}}\big)\to
{\rm H}_0\big(\Tor^{{\cal H}}_1(({\cal H}/{\cal I})^{\otimes\bullet +1},\ad({\cal H}))\big) \to \cdots
\\
\cdots\to \Tor^{{\cal H}}_1(k, \ad({\cal H})) \to {\rm HH}_1\big({\cal A}\mid {\cal H}^{\co{\cal H}/{\cal I}}\big) \to 0.
\end{gather*}

\smallskip

The classical result follows from this spectral sequence by the degeneration argument.

Furthermore, inverting arrows in all diagrams def\/ining the cyclic-homological Takeuchi--Galois transform, interchanging
everywhere \emph{left} and \emph{right}, next applying Connes' cyclic duality and f\/inally inverting the resulting isomorphism we obtain

\smallskip

\noindent\textbf{Theorem~3.7.} \emph{Let ${\cal B}\subseteq {\cal H}$ be a~left comodule subalgebra in a~Hopf algebra
${\cal H}$ such that ${\cal H}\rightarrow {\cal H}/{\cal B}^{+}{\cal H}$ is a~${\cal B}$-Galois coextension, and ${\rm
coad}({\cal H})={\cal H}$ be the right-left SAYD module with the right action given by the multiplication of ${\cal H}$,
and the left coadjoint coaction.
Then there exists an isomorphism
\begin{gather*}
\gamma_{n}\colon \ {\rm C}_{n}({\cal B}, {\rm coad}({\cal H}))^{{\cal H}} \longrightarrow {\rm C}_{n}({\cal H}\mid {\cal H}/{\cal B}^{+}{\cal H})
\end{gather*}
of cyclic modules, defined by~\eqref{aux-vp-psi},~\eqref{aux-vprime-psiprime}.}

\smallskip

Since the domain of the isomorphism depends only on the left comodule subalgebra ${\cal B}$ of a~Hopf algebra ${\cal
H}$, and the codomain depends on the right ${\cal H}$-module quotient coalgebra ${\cal H}/{\cal B}^{+}{\cal H}$ of
${\cal H}$, it can be viewed as a~\emph{cyclic-homological dual Takeuchi--Galois transform} accompanying the
Takeuchi--Galois transform ${\cal B}\mapsto {\cal H}/{\cal B}^{+}{\cal H}$.

Surprisingly, this dual picture is even more interesting than the original construction directly motivated by classical
geometry.
The point is that it makes sense in classical geometry as well, but then is related to quite nontrivial phenomena
connecting geometry and representation theory.
In Subsection 3.5.2 we discuss relations of the dual picture to the \emph{Chern character}, the \emph{Frobenius
reciprocity} and \emph{inertia phenomena} in the \emph{local Langlands program}, the \emph{Chen--Ruan--Brylinski--Nistor
orbifold cohomology} and the \emph{Clifford theory}.

\section{Preliminaries}

In this section we recall the material we will use in the sequel.
In the f\/irst subsection we recall the coalgebra-Galois extensions, as well as the algebra-Galois coextensions.
In the second subsection we discuss the Pontryagin like duality between discrete and linearly compact vector spaces.
Upon recalling the Takeuchi--Galois correspondence between left comodule subalgebras and right module quotient
coalgebras of a~Hopf algebra, we show Pontryagin self-duality of this correspondence.
In the third subsection we recall the relative cyclic homology of algebra extensions.
Finally, the Hopf-cyclic homology, with coef\/f\/icients, of coalgebras is recalled in the fourth subsection.

Throughout the paper, all algebras, coalgebras, and Hopf algebras are over a~f\/ield~$k$, and similarly all unadorned
tensor product symbols $\otimes$ are also over~$k$.
Suppressing the summation, we will write $\Delta(c)=c\ps{1}\otimes c\ps{2}$ for a~comultiplication,
$\rho(v)=v\ns{0}\otimes v\ns{1}$ for a~right coaction, and $\lambda(v)=v\ns{-1}\otimes v\ns{0}$ for a~left coaction.

\subsection{Coalgebra-Galois extensions and algebra-Galois coextensions}

In this subsection we recall the def\/inition and basic properties of the coalgebra-Galois extensions, and the
algebra-Galois coextensions from~\cite{BrzeHaja99,BrzeHaja09}.
Fix the base f\/ield~$k$.

We begin with the coalgebra-Galois extensions.
Let ${\cal C}$ be a~coalgebra coaugmented by the choice of a~group-like element $e\in {\cal C}$.
We say that an algebra ${\cal A}$ is a~\emph{${\cal C}$-algebra} if there is given a~right entwining $\psi\colon  {\cal
C}\otimes {\cal A}\rightarrow {\cal A}\otimes {\cal C}$.
Then ${\cal A}$ becomes a~right ${\cal C}$-comodule via $\rho\colon {\cal A} \to {\cal A} \otimes {\cal C}$ def\/ined as
$\rho(a) =a\ns{0} \otimes a\ns{1}:=\psi(e\otimes a)$.
Then the subspace ${\cal A}^{\co {\cal C}}$ of \emph{coaction invariants} of ${\cal A}$ def\/ined as
\begin{gather*}
{\cal A}^{\co {\cal C}}:=\{b \in {\cal A} \mid \psi(e\otimes b)=b\otimes e\}
\end{gather*}
is a~subalgebra of ${\cal A}$.
Note that ${\cal A}^{\co {\cal C}}$ can be def\/ined as a~limit of a~diagram in vector spaces
\begin{gather*}
{\cal A}^{\co{\cal C}}={\rm Eq}({\cal A}\rightrightarrows {\cal A}\otimes {\cal C})
\end{gather*}
consisting of a~pair of maps $a\mapsto a\ns{0} \otimes a\ns{1}$ and $a\mapsto a\otimes e$.

An algebra extension ${\cal B}\to{\cal A}$ is called a~\emph{${\cal C}$-extension} if ${\cal A}$ is a~${\cal C}$-algebra
and ${\cal B} ={\cal A}^{\co {\cal C}}$.
Finally, a~${\cal C}$-extension ${\cal B}\to{\cal A}$ is said to be \emph{Galois} if the left ${\cal A}$-linear right
${\cal C}$-colinear map
\begin{gather}
\label{aux-canonical-map}
\operatorname{can}\colon \ {\cal A} \otimes_{\cal B} {\cal A} \longrightarrow {\cal A} \otimes {\cal C},
\qquad
a\otimes_{\cal B} a' \mapsto aa'\ns{0} \otimes a'\ns{1}
\end{gather}
is bijective.
The map~\eqref{aux-canonical-map} is called the \emph{canonical map} of the ${\cal C}$-extension.

A \emph{standard} example of a~${\cal C}$-algebra ${\cal A}$ comes from any coaugmented right ${\cal H}$-module
coalgebra~${\cal C}$ and any right ${\cal H}$-comodule algebra ${\cal A}$ over a~Hopf algebra ${\cal H}$.
If we denote by $a\mapsto a_{[0]}\otimes a_{[1]}$ the right ${\cal H}$-coaction on ${\cal A}$ the corresponding
entwining reads as $\psi(c\otimes a)=a_{[0]}\otimes c\cdot a_{[1]}$.

An interesting, due to the geometric examples it covers, class of standard coalgebra-Galois extensions consists of so
called \emph{quotient coalgebra-Galois extensions} def\/ined as follows.
In this setting, one lets~${\cal H}$ to be a~Hopf algebra, a~coaugmented right ${\cal H}$-module coalgebra ${\cal C}$
a~quotient~${\cal H}/{\cal I}$ of~${\cal H}$ by a~coideal right ideal ${\cal I} \subseteq {\cal H}$, and ${\cal A}$
a~right ${\cal H}$-comodule algebra via $\widetilde{\rho}\colon {\cal A} \longrightarrow {\cal A} \otimes {\cal H}$.
Then the composition
\begin{gather}
\label{rho}
\xymatrix{
\rho\colon \ {\cal A} \ar[r]^-{\widetilde{\rho}}&{\cal A} \otimes {\cal H} \ar[r]^-{{\cal A} \otimes\pi}&{\cal A} \otimes {\cal H}/{\cal I}}
\end{gather}
expresses the standard right ${\cal H}/{\cal I}$-coaction on ${\cal A}$.
Therefore the subalgebra of the coaction invariants can be expressed diagrammatically in vector spaces as the equalizer
\begin{gather}
\label{diageq}
{\cal A}^{\co{\cal H}/{\cal I}}={\rm Eq}({\cal A}\rightrightarrows {\cal A}\otimes{\cal H}/{\cal I}),
\end{gather}
where one arrow is~$\rho$ def\/ined by~\eqref{rho}, and the other is the composition
\begin{gather*}
\xymatrix{{\cal A} \ar[r]^-{\cong}&{\cal A} \otimes k\ar[r]^-{{\cal A} \otimes\eta}&{\cal A} \otimes {\cal
H}\ar[r]^-{{\cal A} \otimes\pi} &  {\cal A} \otimes {\cal H}/{\cal I}.}
\end{gather*}
Finally, a~Galois ${\cal H}/{\cal I}$-extension ${\cal A}^{\co {\cal H}/{\cal I}} \to {\cal A}$ is called
a~\emph{quotient coalgebra-Galois extension}.
If~${\cal I}$ is a~Hopf ideal, we call such an extension \emph{quotient Hopf--Galois extension}.

The quantum instanton bundle of~\cite{BoneCiccDabrTarl04,BoneCiccTarl02} is a~quotient coalgebra-Galois extension, and
it uses the full generality of the quotient coalgebra-Galois extensions as ${\cal A} \neq {\cal H}$ and ${\cal I} \neq 0$.

For ${\cal I} =0$, the quotient coalgebra-Galois extensions recover the Hopf--Galois extensions, and in case of ${\cal
A} ={\cal H}$ they are called \emph{homogeneous coalgebra-Galois extensions}.
In particular, viewing ${\cal H}$ as a~right ${\cal H}$-comodule algebra via its comultiplication, one obtains the
homogeneous ${\cal H}/{\cal I}$-Galois extensions
\begin{gather*}
\xymatrix{
\rho\colon \ {\cal H} \ar[r]^{\Delta}& {\cal H} \otimes {\cal H} \ar[r]^{{\cal H} \otimes \pi\ \ \ \ \ \ \ \ \ \ \ \ \ \ } &
{\cal H}\otimes {\cal H}/{\cal I}, \qquad  h \mapsto h\ps{1} \otimes \wbar{h\ps{2}}.
}
\end{gather*}
In the context of the faithfully f\/lat homogeneous extensions, the diagrammatical def\/inition of invariants of the
coaction~\eqref{diageq} plays an important role in the context of Pontryagin self-duality of the Takeuchi--Galois
correspondence considered in the next subsection.

Such Galois extensions correspond to the quantum homogeneous spaces considering ${\cal H}={\cal O}(G)$ the algebra of
functions on a~quantum group $G={\rm Spec}({\cal H})$, and ${\cal B}={\cal O}(X)$ the algebra of functions on a~quantum
space $X={\rm Spec}({\cal B})$.
Then the~$G$-action on~$X$ is encoded by the ${\cal H}$-coaction~\cite[Section~1]{DijkKoor94}.
There is a~celebrated example of this construction due to Podle\'s~\cite{Brze96,Podl87} which is a~quantum spherical
f\/ibration $SU_q(2)\rightarrow S^2_{q,\mu, \nu}$, see~\cite{Brze97}, a~quantum deformation of the classical Hopf
f\/ibration of a~3-sphere over a~2-sphere into circles.

We now recall the algebra-Galois coextensions.
Let ${\cal B}$ be an algebra augmented by the choice of a~character $y\colon  {\cal B}\rightarrow k$.
We say that a~coalgebra ${\cal D}$ is a~\emph{${\cal B}$-coalgebra} if there is given a~left entwining $\varphi\colon {\cal
B}\otimes {\cal D}\rightarrow {\cal D}\otimes {\cal B}$.
Then ${\cal D}$ becomes a~left ${\cal B}$-module via $\ell\colon {\cal B} \otimes {\cal D} \to {\cal D}$ def\/ined as
$\ell=({\cal D}\otimes y)\circ \varphi$.
Then the quotient space ${\cal D}/{\cal B}^{+}{\cal D}$ of \emph{action coinvariants} of ${\cal D}$ def\/ined~by
\begin{gather*}
{\cal B}^{+}:=\{b \in {\cal B} \mid y(b)=0\}
\end{gather*}
is a~quotient coalgebra of ${\cal D}$.
Note that ${\cal D}/{\cal B}^{+}{\cal D}$ can be def\/ined as a~colimit of a~diagram in vector spaces
\begin{gather*}
{\cal D}/{\cal B}^{+}{\cal D}={\rm Coeq}({\cal D}\leftleftarrows {\cal B}\otimes {\cal D})
\end{gather*}
consisting of a~pair of maps $\ell\colon b\otimes d\mapsto b\cdot d$ and $b\otimes d\mapsto y(b)d$.

A coalgebra coextension ${\cal D}\to{\cal C}$ is called a~\emph{${\cal B}$-coextension} if ${\cal D}$ is a~${\cal
B}$-coalgebra and ${\cal C} ={\cal D}/{\cal B}^{+}{\cal D}$.
Finally, a~${\cal B}$-coextension ${\cal D}\to{\cal C}$ is said to be \emph{Galois} if the left ${\cal B}$-module, right
${\cal D}$-comodule map
\begin{gather}
\label{aux-cocan}
\operatorname{cocan}\colon \ {\cal B}\otimes {\cal D}\longrightarrow {\cal D}\Box_{\cal C}{\cal D},
\qquad
b\otimes d\mapsto b\cdot d\ps{1}\otimes d\ps{2}
\end{gather}
is bijective.
The map~\eqref{aux-cocan} is called the \emph{cocanonical map} of the coalgebra ${\cal B}$-coextension.

It is evident that the notion of \emph{algebra-Galois coalgebra coextension} dualizes the notion of
\emph{coalgebra-Galois algebra extension} by formal inverting all arrows in all diagrams and interchanging \emph{left}
and \emph{right}.
In the Hopf--Galois setting the construction dualizes the Hopf--Galois extensions~\cite{Schn93}.

\subsection{Formal Pontryagin duality and Takeuchi--Galois correspondence}\label{subsect-Pontryagin-dual}

In this subsection we will f\/irst summarize the basic properties of the dualization functor on vector spaces
from~\cite{BergmanHausknecht-book}.
We will then recall a~one-to-one correspondence between the coideal subalgebras and quotient coalgebras, known as the
Takeuchi--Galois correspondence~\cite{Take72}.

From the point of view of linear topological vector spaces with continuous linear mappings as morphisms, the dualization
functor def\/ines an equivalence between the opposite category of discrete vector spaces and the category of linearly
compact vector spaces~\cite[Proposition~24.8]{BergmanHausknecht-book}.
Moreover, it transforms naturally the algebraic tensor product of discrete vector spaces into a~completed one of
linearly compact spaces~\cite[Corollary~24.25]{BergmanHausknecht-book}.
In other words, dualization is a~strong monoidal functor.
Therefore one can regard on equal footing all structures def\/ined by diagrams in vector spaces together with their dual
counterparts obtained by reversing all arrows in all necessary diagrams.
This regards linear subspaces and quotient spaces, (co)algebras and Hopf algebras, as well as their bi(co)modules or one
sided (co)modules, one-sided and two-sided (co)ideals, and their (co)tensor products.
Hence, having a~diagrammatical proof of a~theorem in a~symmetric monoidal category of linear (topological) spaces, one
has automatically a~dual theorem after an appropriate dualization of the structures.
It is well known that the notion of Hopf algebra is self-dual.
In particular, it transforms the Hopf group-algebra of a~discrete group into a~linearly compact topological Hopf algebra
of functions on that group, customarily regarded as a~group algebra of a~dual compact quantum group.
Therefore we will call this duality \emph{Pontryagin} to separate it from cyclic duality.
We will show that the notion of SAYD module is Pontryagin self-dual up to the interchange of \emph{left} and
\emph{right}.
According to~\cite{Brze97,BrzeHaja99} the notion of coalgebra-Galois extension of algebras is Pontryagin dual to the
notion of algebra-Galois coextension of coalgebras.
Both dualities play a~role in the present paper.

For any Hopf algebra ${\cal H}$ with multiplication $\mu\colon {\cal H} \otimes{\cal H}\longrightarrow {\cal H}$, unit $\eta\colon
k\longrightarrow {\cal H}$, comultiplication $\Delta\colon {\cal H}\longrightarrow{\cal H}\otimes{\cal H}$, counit
$\varepsilon\colon {\cal H}\longrightarrow k$, and antipode $S\colon {\cal H}\longrightarrow{\cal H}$, Takeuchi introduces
in~\cite{Take72} the one-to-one (Galois) correspondence
\begin{gather*}
\{{\cal B}\subseteq {\cal H} \mid {\cal B}\ {\rm is\ a\ left\ coideal\ subalgebra}, {\cal H}\ {\rm is\ faithfully\ f\/lat\
over}\ {\cal B} \}
\\
\hspace*{50mm}\uparrow\downarrow
\\
\{{\cal I}\subseteq {\cal H} \mid {\cal I}\ {\rm is\ a\ coideal\ right\ ideal}, {\cal H}\ {\rm is\ faithfully\ cof\/lat\
over}\ {\cal H}/{\cal I} \},
\end{gather*}
under which ${\cal B}\mapsto {\cal B}^{+}{\cal H}$ and ${\cal I}\mapsto {\cal H}^{\co {\cal H}/{\cal I}}$.

A crucial observation for our purpose is that the above correspondence can be written categorically as an equivalence
between the category of left ${\cal H}$-comodule f\/lat extensions $i\colon {\cal B}\to{\cal H}$ and the category of right
${\cal H}$-module cof\/lat coextensions $\pi\colon {\cal H}\to {\cal C}$, which is given~by
\begin{gather*}
{\cal B}\mapsto {\rm Coeq}({\cal H}\leftleftarrows {\cal B}\otimes{\cal H}),
\qquad
{\cal C}\mapsto {\rm Eq}({\cal H}\rightrightarrows {\cal H}\otimes{\cal C}),
\end{gather*}
where the parallel pair of left arrows in the coequalizer and the parallel pair of right arrows in the equalizer read as composites
\begin{gather*}
\xymatrix@R-1pc{{\cal H} & {\cal H}\otimes{\cal H}\ar[l]_-{\mu}& {\cal B}\otimes{\cal H}\ar[l]_-{\iota\otimes {\cal H}},
\\
{\cal H} & k\otimes {\cal H}\ar[l]_-{\cong} &{\cal B}\otimes{\cal H}\ar[l]_-{\varepsilon\iota\otimes {\cal H}},}\qquad
 \xymatrix@R-1pc{{\cal H}\ar[r]^-{\Delta} & {\cal H}\otimes{\cal H} \ar[r]^-{{\cal H}\otimes\pi} &{\cal H}\otimes{\cal C},
\\
{\cal H} \ar[r]^-{\cong} & {\cal H}\otimes k \ar[r]^-{{\cal H}\otimes \pi\eta} &{\cal H}\otimes{\cal C},}
\end{gather*}
respectively.
It is thus evident that this correspondence is Pontryagin self-dual up to an interchange of \emph{left} and \emph{right}
in dual structures.
We will say that ${\cal B}$ and ${\cal C}$ as above are \emph{Takeuchi--Galois transforms} of each other.

\subsection{Relative cyclic homology of (co)algebra (co)extensions}

We start with a~quick detour on the relative cyclic homology of algebra extensions.
For any ${\cal B}$-bimodule ${\cal M}$ we let
\begin{gather*}
[{\cal M}]_{\cal B}:={\cal M}/[{\cal M},{\cal B}],
\end{gather*}
where $[{\cal M},{\cal B}]$ is the subspace of ${\cal M}$ generated by all commutators $[m,b]:=m\cdot b-b\cdot m$.
We will denote by $[m]_{\cal B}$ the class of $m\in {\cal M}$ in $[{\cal M}]_{\cal B}$.

The relative cyclic homology of an algebra extension ${\cal B}\rightarrow {\cal A}$ is computed by the cyclic object
\begin{gather}
\label{aux-rel-cyclic-alg}
{\rm C}_n({\cal A}\mid{\cal B}):=[{\cal A}^{\otimes_{\cal B} n+1}]_{\cal B},
\end{gather}
equipped, for all $n\geq 0$, with the morphisms
\begin{gather*}
d_i\colon \ {\rm C}_n({\cal A}\mid{\cal B}) \longrightarrow {\rm C}_{n-1}({\cal A}\mid{\cal B}),
\qquad
0 \leq j \leq n,
\\
d_i\big[a^0 \otimes_{\cal B} a^1 \otimes_{\cal B} \dots \otimes_{\cal B} a^n\big]_{\cal B}\\
\qquad{}
=\begin{cases}
\big[a^0 a^1 \otimes_{\cal B} \dots \otimes_{\cal B} a^n \big]_{\cal B}, & i =0,
\\
\big[a^0 \otimes_{\cal B} a^1 \otimes_{\cal B} \dots \otimes_{\cal B} a^i a^{i+1} \otimes_{\cal B} \dots \otimes_{\cal B} a^n \big]_{\cal B},& 1\leq i\leq n-1,
\\
\big[a^n a^0 \otimes_{\cal B} a^1 \otimes_{\cal B} \dots \otimes_{\cal B} a^{n-1}\big]_{\cal B}, & i =n,
\end{cases}
\\
s_i\colon \ {\rm C}_n({\cal A}\mid{\cal B}) \longrightarrow {\rm C}_{n+1}({\cal A}\mid{\cal B}),
\qquad
0 \leq j \leq n,
\\
s_i\big[a^0 \otimes_{\cal B} \dots \otimes_{\cal B} a^n\big]_{\cal B} =\big[a^0 \otimes_{\cal B} a^1
\otimes_{\cal B} \dots \otimes_{\cal B} a^i \otimes_{\cal B} 1 \otimes_{\cal B} a^{i+1} \otimes_{\cal B} \dots
\otimes_{\cal B} a^n \big]_{\cal B}
\end{gather*}
and
\begin{gather*}
t_n\colon \ {\rm C}_n({\cal A}\mid{\cal B}) \longrightarrow {\rm C}_n({\cal A}\mid{\cal B}),
\qquad
t_n\big[a^0 \otimes_{\cal B} \dots \otimes_{\cal B} a^n\big]_{\cal B} =\big[a^n \otimes_{\cal B} a^0
\otimes_{\cal B} \dots \otimes_{\cal B} a^{n-1}\big]_{\cal B}.
\end{gather*}
Hochschild homology of the complex~\eqref{aux-rel-cyclic-alg} is called the Hochschild homology of the extension, and is
denoted by ${\rm HH}_\bullet({\cal A}\mid{\cal B})$.
Similarly, cyclic (resp.\
periodic cyclic, negative cyclic) homology of the cyclic object~\eqref{aux-rel-cyclic-alg} is called the relative cyclic
(resp.\ Hochschild, periodic cyclic, negative cyclic) homology of the extension ${\cal B} \rightarrow{\cal A}$, and it is
denoted by ${\rm HC}_\bullet({\cal A}\mid{\cal B})$ (resp.\
${\rm HH}_\bullet({\cal A}\mid{\cal B})$, ${\rm HP}_\bullet({\cal A}\mid{\cal B})$, ${\rm HN}_\bullet({\cal A}\mid{\cal B})$).

The natural mapping
\begin{gather*}
{\rm C}_n({\cal A})\longrightarrow {\rm C}_n({\cal A}\mid{\cal B}),
\qquad
a^0\otimes\dots\otimes a^n\mapsto \big[a^0\otimes_{\cal B} \dots \otimes_{\cal B} a^n\big]_{{\cal B}},
\qquad
n\geq 0,
\end{gather*}
induces a~map ${\rm HH}_\bullet({\cal A})\longrightarrow {\rm HH}_\bullet({\cal A}\mid{\cal B})$ of Hochschild homology
groups, which further induces a~mapping ${\rm HC}_\bullet({\cal A})\longrightarrow {\rm HC}_\bullet({\cal A}\mid{\cal
B})$ of cyclic homology groups.
This is the computational motivation behind the relative cyclic homology~\cite{Kadi92}, in which it is proved to be an
isomorphism when ${\cal B}\subseteq {\cal A}$ is a~separable subalgebra (semisimple in characteristic zero).
The relative cyclic homology then appeared in~\cite{Scha92}, where the relative cyclic homology ${\rm HC}_\bullet(kG\mid kN)$
of group algebras, associated to normal subgroups, was computed, and is applied to extend Eckmann's
result~\cite{Eckm86} on the Bass conjecture.

Dually, for any ${\cal C}$-bicomodule ${\cal M}$ one def\/ines ${\cal M}^{\cal C}:=\{m\in {\cal M}\mid m\ns{0}\otimes
m\ns{1} =m\ns{0}\otimes m\ns{-1}\}$. Then the relative cyclic cohomology of a~coalgebra coextension ${\cal
D}\rightarrow {\cal C}$ is computed by the cocyclic object
\begin{gather}
\label{aux-rel-cyclic-coalg}
{\rm C}^n({\cal D}\mid{\cal C}):=\big({\cal D}^{\Box_{\cal C} n+1}\big)^{\cal C},
\qquad
n\geq 0,
\end{gather}
with the structure maps
\begin{gather*}
\delta_i\colon \ {\rm C}^n({\cal D}\mid{\cal C}) \longrightarrow {\rm C}^{n+1}({\cal D}\mid{\cal C}),
\qquad
0 \leq j \leq n+1,
\\
\delta_i\big(d^0\otimes d^1\otimes \cdots\otimes d^n\big)
=\begin{cases}
d^0\ps{1}\otimes d^0\ps{2}\otimes d^1\otimes \cdots\otimes d^n, & i =0,
\\
d^0\otimes \cdots \otimes  d^i\ps{1}\otimes d^i\ps{2}\otimes \cdots \otimes d^n, & 1 \leq i\leq n,
\\
d^0\ps{2}\otimes d^1\otimes \cdots\otimes d^n\otimes d^0\ps{1}, & i =n+1,
\end{cases}
\\
\sigma_j\colon \ {\rm C}^n({\cal D}\mid{\cal C}) \longrightarrow {\rm C}^{n-1}({\cal D}\mid{\cal C}),
\qquad
0 \leq j \leq n-1,
\\
\sigma_j\big(d^0\otimes d^1\otimes \cdots\otimes d^n\big) =d^0\otimes d^1\otimes \cdots \otimes
d^{j}\varepsilon\big(d^{j+1}\big)\otimes \cdots \otimes d^n
\end{gather*}
and
\begin{gather*}
\tau_n\colon \ {\rm C}^n({\cal D}\mid{\cal C}) \longrightarrow {\rm C}^n({\cal D}\mid{\cal C}),
\qquad
\tau_n\big(d^0\otimes d^1\otimes \cdots\otimes d^n\big) =d^1\otimes \cdots\otimes d^n\otimes d^0.
\end{gather*}
Then the cyclic duality~\cite{Conn83,KhalRang05} yields the cyclic module structure on the collection of vector spaces
\begin{gather}
\label{cyc-mod-coext}
{\rm C}_n({\cal D}\mid{\cal C}):=\big({\cal D}^{\Box_{\cal C} n+1}\big)^{\cal C}
\end{gather}
given by faces
\begin{gather*}
d_i\colon \ {\rm C}_n({\cal D}\mid{\cal C}) \longrightarrow {\rm C}_{n-1}({\cal D}\mid{\cal C}),
\qquad
0 \leq i \leq n,
\\
d_i\big(d^0\otimes d^1\otimes \cdots\otimes d^n\big) =d^0\otimes d^1\otimes \cdots \otimes \varepsilon\big(d^i\big)d^{i+1}\otimes\cdots \otimes d^n,
\qquad
0 \leq i \leq n-1,
\\
d_n\big(d^0\otimes d^1\otimes \cdots\otimes d^n\big) =d^0\otimes d^1\otimes \cdots \otimes d^{n-1}\varepsilon\big(d^n\big),
\end{gather*}
degeneracies
\begin{gather*}
s_i\colon \ {\rm C}_n({\cal D}\mid{\cal C}) \longrightarrow {\rm C}_{n+1}({\cal D}\mid{\cal C}),
\qquad
0 \leq i \leq n,
\\
s_i\big(d^0\otimes d^1\otimes \cdots\otimes d^n\big) =d^0\otimes d^1\otimes \cdots \otimes\Delta\big(d^i\big)\otimes \cdots \otimes
d^n,
\end{gather*}
and the cyclic operator
\begin{gather*}
t_n\colon \ {\rm C}_n({\cal D}\mid{\cal C}) \longrightarrow {\rm C}_n({\cal D}\mid{\cal C}),
\qquad
t_n\big(d^0\otimes d^1\otimes \cdots\otimes d^n\big) =d^n\otimes d^0\otimes \cdots\otimes d^{n-1}.
\end{gather*}
Hochschild homology of the complex~\eqref{cyc-mod-coext} is called the Hochschild homology of the coextension, and is
denoted by ${\rm HH}_\bullet({\cal D}\mid{\cal C})$.
Similarly, cyclic (resp.\
periodic cyclic, negative cyclic) homology of the cyclic object~\eqref{aux-rel-cyclic-coalg} is called the relative
cyclic (resp.\
Hochschild, periodic cyclic, negative cyclic) homology of the coextension ${\cal D} \rightarrow{\cal C}$, and it is
denoted by ${\rm HC}_\bullet({\cal D}\mid{\cal C})$ (resp.\
${\rm HH}_\bullet({\cal D}\mid{\cal C})$, ${\rm HP}_\bullet({\cal D}\mid{\cal C})$, ${\rm HN}_\bullet({\cal D}\mid{\cal C})$).

\subsection[Hopf-cyclic homology of ${\cal H}$-module coalgebras and ${\cal H}$-comodule algebras]{Hopf-cyclic homology
of $\boldsymbol{\cal H}$-module coalgebras and $\boldsymbol{\cal H}$-comodule algebras}

In this subsection we recall the relative Hopf-cyclic homology with coef\/f\/icients, for coalgebras, using the cyclic
duality principle~\cite{Conn83,KhalRang05}.

Let us f\/irst recall the def\/inition of a~left-right stable anti-Yetter--Drinfeld module over a~Hopf algebra ${\cal H}$
from~\cite{HajaKhalRangSomm04-I}.
Let a~linear space ${\cal M}$ be a~left ${\cal H}$-module by $\ell\colon {\cal H}\otimes{\cal M}\longrightarrow {\cal M}$
given by $\ell(h\otimes m)=h\cdot m$, and a~right ${\cal H}$-comodule via $\rho\colon {\cal M} \longrightarrow {\cal M}
\otimes {\cal H}$ given by $\rho(m)=m\ns{0} \otimes m\ns{1}$.
Then ${\cal M}$ is called a~left-right anti-Yetter--Drinfeld module (AYD module) over ${\cal H}$ if the ${\cal
H}$-action and ${\cal H}$-coaction are compatible as
\begin{gather*}
(h \cdot m)\ns{0}\otimes (h \cdot m)\ns{1} =h\ps{2} \cdot m\ns{0} \otimes h\ps{3} m\ns{1} S(h\ps{1}).
\end{gather*}
A~left-right AYD module ${\cal M}$ is called \emph{stable} (and then is abbreviated as SAYD module) if
\begin{gather*}
m\ns{1}\cdot m\ns{0} =m.
\end{gather*}
On the other hand, a~right-left SAYD module structure is def\/ined in terms of a~right ${\cal H}$-module structure via
$r\colon {\cal M}\otimes{\cal H}\longrightarrow {\cal M}$ given by $r(m\otimes h)=m\cdot h$, and a~left ${\cal H}$-comodule
structure via $\lambda\colon {\cal M} \longrightarrow {\cal H} \otimes {\cal M}$ given by $\lambda(m)=m\ns{-1} \otimes
m\ns{0}$ such that
\begin{gather*}
(m \cdot h)\ns{-1}\otimes (m \cdot h)\ns{0} =S(h\ps{3}) \cdot m\ns{-1} h\ps{1} \otimes m\ns{0} h\ps{2},
\qquad
m\ns{0}\cdot m\ns{-1} =m.
\end{gather*}
Note that these conditions are Pontryagin dual to each other.
Indeed, the left-right SAYD module compatibility reads as commutativity of the following diagrams in the category $V\
ect$ of vector spaces ($\sigma$ being the transposition of the corresponding tensorands)
\begin{gather*}
\xymatrix{{\cal H}\otimes {\cal H}\otimes {\cal H}\otimes {\cal M}\otimes {\cal H}\ar[rr] && {\cal H}\otimes {\cal
M}\otimes {\cal H}\otimes {\cal H}\otimes {\cal H}\ar[d]^-{\ell \otimes \mu_{3}}
\\
{\cal H}\otimes {\cal M}\ar[u]^-{\Delta^3\otimes \rho}\ar[r]^-{\ell} & {\cal M} \ar[r]^-{\rho}&{\cal M}\otimes {\cal H},}
\end{gather*}
where the upper horizontal arrow admits two decompositions
\begin{gather*}
\xymatrix{& {\cal H}\otimes {\cal H}\otimes {\cal M}\otimes {\cal H}\otimes {\cal H} \ar[rd]^{\  \ \
H\otimes\sigma_{{\cal H},{\cal M}}\otimes {\cal H}\otimes S}&
\\
{\cal H}\otimes {\cal H}\otimes {\cal H}\otimes {\cal M}\otimes {\cal H}\ar[ur]^{\sigma_{{\cal H}, {\cal H}\otimes {\cal
H}\otimes {\cal M}\otimes {\cal H}}\ \ \ \ }\ar[dr]_{S\otimes {\cal H} \otimes \sigma_{{\cal H}, {\cal
M}}\otimes {\cal H} \ \ \ \ \ } && {\cal H}\otimes {\cal M}\otimes {\cal H}\otimes {\cal H}\otimes {\cal H}
\\
&{\cal H}\otimes {\cal H}\otimes {\cal M}\otimes {\cal H}\otimes {\cal H}\ar[ru]_{\ \ \ \ \ \sigma_{{\cal H}\otimes
{\cal H}\otimes {\cal M}\otimes {\cal H}, {\cal H}}},}
\end{gather*}
and
\begin{gather*}
\xymatrix{{\cal M}\otimes {\cal H} \ar[r]^-{\sigma_{{\cal M},{\cal H}}} & {\cal H}\otimes {\cal M}\ar[d]^-{\ell}
\\
{\cal M}\ar[u]^{\rho}\ar[r]^{=} & {\cal M}.}
\end{gather*}
Now we see that reversing the arrows, inverting transpositions~$\sigma$ and next interchanging in pairs~$\Delta$
and~$\mu$, $\ell $ and~$\lambda$,~$\rho$ and~$r$, and f\/inally \emph{left} and \emph{right}, we obtain the right-left
SAYD module compatibility.

Let us now recall the Hopf-cyclic cohomology of a~Hopf-module coalgebra with coef\/f\/icients.
Let ${\cal H}$ be a~Hopf algebra with an invertible antipode, ${\cal C}$ a~right ${\cal H}$-module coalgebra, i.e.,
a~right ${\cal H}$-module such that
\begin{gather*}
\Delta(c\cdot h) =c\ps{1}\cdot h\ps{1}\otimes c\ps{2}\cdot h\ps{2},
\qquad
\varepsilon(c\cdot h) =\varepsilon(c)\varepsilon(h),
\qquad
\forall\,  c\in {\cal C}, \quad \forall\,  h\in{\cal H},
\end{gather*}
and ${\cal M}$ a~left-right SAYD module over ${\cal H}$.
Then for all $n\geq 0$ we def\/ine
\begin{gather}
\label{aux-cyclic-module-C-H-M}
{\rm C}^n({\cal C},{\cal M})_{{\cal H}}:={\cal C}^{\otimes n+1} \otimes_{{\cal H}}{\cal M},
\end{gather}
where on the right hand side the diagonal right ${\cal H}$-module structure on ${\cal C}^{\otimes n+1}$ is used.

{\sloppy The subscript ${\cal H}$ refers to the coinvariants of the diagonal right ${\cal H}$-module structure \mbox{${\cal C}^{\otimes
n+1}\otimes{\cal M}$}, by switching the one on ${\cal M}$ from left to right via the antipode.
The collection of vector spaces ${\rm C}^n({\cal C},{\cal M})_{{\cal H}}$ together with the operators
\begin{gather*}
\delta_i\colon \ {\rm C}^{n-1}({\cal C},{\cal M})_{{\cal H}} \to {\rm C}^n({\cal C},{\cal M})_{{\cal H}},
\qquad
0 \leq i \leq n,
\\
\delta_i\big(\big(c^0 \otimes \cdots \otimes c^{n-1}\big) \otimes_{{\cal H}} m\big) \\
\qquad{} =\big(c^0 \otimes \cdots \otimes c^i\ps{1} \otimes
c^i\ps{2} \otimes \cdots \otimes c^{n-1}\big) \otimes_{{\cal H}} m,
\qquad
0 \leq i \leq n-1,
\\
\delta_n\big(\big(c^0 \otimes \cdots \otimes c^{n-1}\big) \otimes_{{\cal H}} m\big) =\big(c^0\ps{2} \otimes \cdots \otimes c^{n-1}
\otimes c^0\ps{1} \cdot S^{-1}(m\ns{1})\big) \otimes_{{\cal H}} m\ns{0},
\\
\sigma_j\colon \ {\rm C}^{n+1}({\cal C},{\cal M})_{{\cal H}} \to {\rm C}^n({\cal C},{\cal M})_{{\cal H}},
\qquad
0 \leq i \leq n,
\\
\sigma_j\big(\big(c^0 \otimes \cdots \otimes c^{n+1}\big) \otimes_{{\cal H}} m\big) =\big(c^0 \otimes \cdots \otimes
c^{j}\varepsilon\big(c^{j+1}\big) \otimes \cdots \otimes c^{n+1}\big) \otimes_{{\cal H}} m,
\end{gather*}
and
\begin{gather*}
\tau_n\colon \ {\rm C}^{n}({\cal C},{\cal M})_{{\cal H}} \to {\rm C}^n({\cal C},{\cal M})_{{\cal H}},
\\
\tau_n\big(\big(c^0 \otimes \cdots \otimes c^n\big) \otimes_{\mathcal{H}} m\big) =\big(c^1 \otimes \cdots \otimes c^n \otimes c^0 \cdot
S^{-1}(m\ns{1})) \otimes_{\mathcal{H}} m\ns{0}.
\end{gather*}
is a~cocyclic module.
Cyclic cohomology of this cocyclic module is called the Hopf-cyclic cohomology of the ${\cal H}$-module coalgebra ${\cal
C}$ with coef\/f\/icients in the left-right SAYD module ${\cal M}$ over~${\cal H}$, and is denoted by ${\rm
HC}^\bullet({\cal C},{\cal M})_{\cal H}$.

}

Applying the cyclic duality procedure~\cite{Conn83,KhalRang05} we obtain on~\eqref{aux-cyclic-module-C-H-M} the cyclic
module structure given by the faces
\begin{gather*}
d_i\colon \ {\rm C}_{n+1}({\cal C},{\cal M})_{{\cal H}} \to {\rm C}_n({\cal C},{\cal M})_{{\cal H}},
\qquad
0 \leq i \leq n+1,
\\
d_i\big(\big(c^0 \otimes \cdots \otimes c^{n+1}) \otimes_{{\cal H}} m\big) =\big(c^0 \otimes \cdots \otimes \varepsilon\big(c^i\big)c^{i+1}
\otimes \cdots \otimes c^{n+1}\big) \otimes_{{\cal H}} m,
\qquad
0 \leq i \leq n,
\\
d_{n+1}\big(\big(c^0 \otimes \cdots \otimes c^{n+1}\big) \otimes_{{\cal H}} m\big) =\big(\varepsilon\big(c^{n+1}\big) c^0 \otimes \cdots \otimes
c^{n}\big) \otimes_{{\cal H}} m
\end{gather*}
degeneracies
\begin{gather*}
s_i\colon \ {\rm C}_{n-1}({\cal C},{\cal M})_{{\cal H}} \to {\rm C}_n({\cal C},{\cal M})_{{\cal H}},
\qquad
0 \leq i \leq n-1,
\\
s_i\big(\big(c^0 \otimes \cdots \otimes c^{n-1}\big) \otimes_{\mathcal{H}} m\big) =\big(c^0 \otimes \cdots \otimes c^i\ps{1} \otimes
c^i\ps{2} \otimes \cdots \otimes c^{n-1}\big) \otimes_{\mathcal{H}} m,
\end{gather*}
and the cyclic operator
\begin{gather*}
t_n\colon \ {\rm C}_{n}({\cal C},{\cal M})_{{\cal H}} \to {\rm C}_n({\cal C},{\cal M})_{{\cal H}},
\\
t_n\big(\big(c^0 \otimes \cdots \otimes c^n\big) \otimes_{{\cal H}} m\big) =\big(c^n \cdot m\ns{1} \otimes c^0 \otimes \cdots \otimes
c^{n-1}\big) \otimes_{{\cal H}} m\ns{0}.
\end{gather*}
The cyclic homology of this cyclic module is called the cyclic (resp.\
Hochschild, periodic cyclic, negative cyclic) homology of the $\mathcal{H}$-module coalgebra ${\cal C}$, with
coef\/f\/icients in the left-right SAYD module ${\cal M}$ over ${\cal H}$, and is denoted by ${\rm HC}_\bullet({\cal
C},{\cal M})_{\cal H}$ (resp.\
${\rm HH}_\bullet({\cal C},{\cal M})_{\cal H}$, ${\rm HP}_\bullet({\cal C},{\cal M})_{\cal H}$, ${\rm HN}_\bullet({\cal
C},{\cal M})_{\cal H}$).

We conclude this section with the Hopf-cyclic homology of ${\cal H}$-comodule algebras~\cite{HajaKhalRangSomm04-II}.
Let ${\cal M}$ be a~left-right SAYD module over the Hopf algebra ${\cal H}$, and ${\cal B}$ a~left ${\cal H}$-comodule algebra.
Then for $n\geq 0$ we def\/ine ${\rm C}_{n}({\cal B},{\cal M})^{{\cal H}}:={\cal M} \Box_{\cal H} {\cal B}^{\otimes n+1}$,
where on the right hand side the diagonal left ${\cal H}$-comodule structure on ${\cal B}^{\otimes n+1}$ is used.
The superscript ${\cal H}$ refers to the invariants of the diagonal right ${\cal H}$-comodule structure ${\cal
M}\otimes{\cal B}^{\otimes n+1}$, switching the one on ${\cal B}^{\otimes n+1}$ from left to right by means of the
antipode.

The collection of these vector spaces together with operators
\begin{gather*}
d_i\colon \ {\rm C}_n({\cal B},{\cal M})^{{\cal H}} \longrightarrow {\rm C}_{n-1}({\cal B},{\cal M})^{{\cal H}},
\qquad
0 \leq i \leq n,
\\
d_i\big(m\otimes b^0 \otimes \cdots \otimes b^n\big) =m\otimes b^0 \otimes \cdots \otimes b^ib^{i+1}\otimes \cdots \otimes b^n,
\qquad
0 \leq i \leq n-1,
\\
d_n\big(m\otimes b^0 \otimes \cdots \otimes b^n\big) =b^n\ns{-1}\cdot m\otimes b^n\ns{0}b^0\otimes b^1 \otimes \cdots \otimes b^{n-1},
\\
s_j\colon \ {\rm C}_n({\cal B},{\cal M})^{{\cal H}} \longrightarrow {\rm C}_{n+1}({\cal B},{\cal M})^{{\cal H}},
\qquad
0 \leq i \leq n,
\\
s_j\big(m\otimes b^0 \otimes \cdots \otimes b^n\big) =m\otimes b^0 \otimes\cdots \otimes b^j\otimes 1\otimes\cdots \otimes b^n,
\end{gather*}
and
\begin{gather*}
t_n\colon \ {\rm C}_n({\cal B},{\cal M})^{{\cal H}} \longrightarrow {\rm C}_n({\cal B},{\cal M})^{{\cal H}},
\\
t_n\big(m\otimes b^0 \otimes \cdots \otimes b^n\big) =b^n\ns{-1}\cdot m\otimes b^n\ns{0}\otimes b^0\otimes \cdots \otimes b^{n-1}
\end{gather*}
is a~cyclic module.
Finally, the cyclic homology of this cyclic module is called the cyclic (resp.\
Hochschild, periodic cyclic, negative cyclic) homology of the $\mathcal{H}$-comodule algebra ${\cal C}$, with
coef\/f\/icients in the left-right SAYD module ${\cal M}$ over ${\cal H}$, and is denoted by ${\rm HC}_\bullet({\cal
C},{\cal M})^{\cal H}$ (resp.\
${\rm HH}_\bullet({\cal C},{\cal M})^{\cal H}$, ${\rm HP}_\bullet({\cal C},{\cal M})^{\cal H}$, ${\rm HN}_\bullet({\cal
C},{\cal M})^{\cal H}$).

\section{Relative cyclic homology as Hopf-cyclic homology\\ with coef\/f\/icients}

In this section we achieve our main result identifying the relative cyclic homology of a~homogeneous ${\cal C}$-Galois
extension ${\cal B}\subseteq{\cal H}$ with the Hopf-cyclic homology, with coef\/f\/icients, of the right ${\cal H}$-module
coalgebra~${\cal C}$.
Moreover, in view of the Pontryagin duality of Subsection~\ref{subsect-Pontryagin-dual}, we obtain an identif\/ication of
the relative cyclic homology of a~${\cal B}$-Galois coextension ${\cal H}\to {\cal C}$ with the Hopf-cyclic homology,
with (the Pontryagin dual) coef\/f\/icients, of the ${\cal H}$-comodule algebra ${\cal B}$.
We shall conclude the section developing spectral sequences to shed further light on the relative homology groups.

\subsection{The isomorphism of cyclic objects}\label{subsect-iso-cyclic-obj}

In this subsection we will construct an explicit isomorphism from the relative homology complex of a~homogeneous ${\cal
H}/{\cal I}$-Galois extension ${\cal B}:={\cal H}^{\co {\cal H}/{\cal I}}\to {\cal H}$ to the Hopf-cyclic homology
complex of the ${\cal H}$-module coalgebra ${\cal H}/{\cal I}$.

Let ${\cal H}$ be a~Hopf algebra and ${\cal I} \subseteq {\cal H}$ a~coideal right ideal of ${\cal H}$.
Then ${\cal H}/{\cal I}$ becomes a~coaugmented quotient coalgebra in a~canonical way,
\begin{gather*}
\Delta(\wbar{h}):=\overline{h\ps{1}}\otimes \overline{h\ps{2}},
\qquad
\varepsilon(\wbar{h}):=\varepsilon(h),
\qquad
\forall\,  h\in{\cal H},
\end{gather*}
where $\wbar{h}:=h+{\cal I}$, and the canonical coaugmentation of ${\cal H}/{\cal I}$ is given by the group-like
$\overline{1}$.

Let ${\cal B} \subseteq {\cal H}$ be a~homogeneous ${\cal H}/{\cal I}$-Galois extension given by the canonical right
${\cal H}/{\cal I}$-coaction
\begin{gather*}
\xymatrix{{\cal H}\ar[r]^{\Delta} & {\cal H}\otimes {\cal H}\ar[r] &{\cal H} \otimes ({\cal H}/{\cal I}),
\qquad
h \mapsto h\ps{1} \otimes \wbar{h\ps{2}}}
\end{gather*}
on ${\cal H}$.
In view of the def\/inition ${\cal B}:={\cal H}^{\co {\cal H}/{\cal I}}$, which reads as
\begin{gather}
\label{coinv}
b\ps{1}\otimes \overline{b\ps{2}} =b\otimes \overline{1},
\end{gather}
applying the counit to the left tensorands we deduce $\overline{b}=\varepsilon(b)\overline{1}$.
We will use also the following iterated form of~\eqref{coinv}:
\begin{gather}
\label{coinviter}
b\ps{1}\otimes \dots\otimes\overline{b\ps{n+2}} =b\ps{1}\otimes \dots\otimes b\ps{n+1}\otimes\overline{1}.
\end{gather}

The three following facts are crucial to our purposes.
First of all, ${\cal B}$ is the left coideal of~${\cal H}$~\cite{BrzeHaja09}, i.e.,
\begin{gather*}
\Delta({\cal B})\subset {\cal H}\otimes {\cal B}.
\end{gather*}
Secondly, the Galois condition f\/ixes the coideal right ideal ${\cal I}\subseteq {\cal H}$ completely as
follows~\cite[Theorem~2.6]{BrzeHaja09}.

\begin{Theorem}
Let ${\cal B}\subseteq {\cal H}$ be a~homogeneous ${\cal H}/{\cal I}$-extension.
Then this extension is Galois if and only if ${\cal I} ={\cal B}^{+}{\cal H}$, where ${\cal B}^{+}:={\cal B} \cap \operatorname{Ker}
\varepsilon$.
\end{Theorem}

Finally, there is an explicit formula for the translation map, and hence for the inverse to the canonical map,~\cite[Corollary~2.8]{BrzeHaja09}.

\begin{Corollary}
Let ${\cal B}\subseteq {\cal H}$ be a~coalgebra-Galois ${\cal H}/{\cal I}$-extension as above.
Then the translation map $\tau:=\operatorname{can}^{-1}(1 \otimes -)$ is given~by
\begin{gather}
\label{trans}
\tau(\overline{h}) =S(h\ps{1}) \otimes_{B} h\ps{2}.
\end{gather}
\end{Corollary}

Moreover, it implies that the right hand side of~\eqref{trans} is independent of the choice of the representative~$h$ of
the class $\overline{h}$.

The canonical map
\begin{gather*}
{\cal H}\otimes_{{\cal B}}{\cal H}\rightarrow {\cal H}\otimes{\cal H}/{\cal I},
\qquad
h\otimes_{{\cal B}}h' \mapsto hh'\ps{1}\otimes \overline{h'\ps{2}}
\end{gather*}
and its inverse
\begin{gather*}
{\cal H}\otimes{\cal H}/{\cal I} \rightarrow {\cal H}\otimes_{{\cal B}}{\cal H},
\qquad
h\otimes \overline{h'}\mapsto h S(h'\ps{1})\otimes_{{\cal B}} h'\ps{2}
\end{gather*}
can be inductively extended to ${\cal H}$-bimodule isomorphisms
\begin{gather*}
\operatorname{can}_n\colon \ {\cal M}\otimes_{\cal B}\otimes{\cal H}^{\otimes_{\cal B} n} \longrightarrow {\cal M} \otimes {\cal C}^{\otimes n},
\\
\operatorname{can}_n\big(m \otimes_{\cal B} h^1\otimes_{\cal B} \cdots \otimes_{\cal B} h^n\big)
\\
\qquad{}
=mh^1\ps{1}\dots h^n\ps{1} \otimes \wbar{h^1\ps{2}\cdots h^n\ps{2}} \otimes\cdots\otimes \wbar{h^{n-1}\ps{n}h^n\ps{n}}
\otimes \wbar{h^n\ps{n+1}}.
\end{gather*}
and
\begin{gather*}
\operatorname{can}^{-1}_n\colon \  {\cal M} \otimes {\cal C}^{\otimes n} \longrightarrow {\cal M}\otimes_{\cal B} {\cal H}^{\otimes_{\cal B} n},
\\
\operatorname{can}^{-1}_n\big(m \otimes \wbar{g^1} \otimes\cdots\otimes \wbar{g^n}\big)
\\
\qquad{}
=mS\big(g^1\ps{1}\big) \otimes_{\cal B} g^1\ps{2}S\big(g^2\ps{1}\big) \otimes_{\cal B} \cdots \otimes_{\cal B} g^{n-1}\ps{2}S\big(g^n\ps{1}\big)
\otimes_{\cal B} g^n\ps{2},
\end{gather*}
respectively, for any ${\cal H}$-bimodule ${\cal M}$, see also~\cite[Proposition~3.6]{Kadi14}.

We also recall that ${\cal H}$, equipped with the (left) adjoint ${\cal H}$-action
\begin{gather}
\label{HonM}
h\triangleright h':=h\ps{2}h'S(h\ps{1})
\end{gather}
and the (right) ${\cal H}$-coaction given by the comultiplication, is a~left-right SAYD module over ${\cal H}$, which is
denoted by $\ad({\cal H})$, see for instance~\cite[Example~4.3]{JaraStef06}.
Moreover, this action satisf\/ies
\begin{gather}
\label{HonH}
h\ps{2}\triangleright (h'h\ps{1}) =h\ps{2}\ps{2}h'h\ps{1}S(h\ps{2}\ps{1})
=h\ps{3}h'h\ps{1}S(h\ps{2}) =h\ps{2}h'\varepsilon(h\ps{1}) =hh'.
\end{gather}

All that is used in the following lemma.

\begin{Lemma}
Let ${\cal B}\subseteq {\cal H}$ be a~homogeneous ${\cal H}/{\cal I}$-Galois extension.
Then for any $n \geq 0$ we have the isomorphism of vector spaces $[{\cal H}^{\otimes_{\cal B} n+1}]_{{\cal
B}} \cong ({\cal H}/{\cal I})^{\otimes n+1} \otimes_{\cal H} \ad({\cal H})$ implemented~by
\begin{gather}
\psi_n\colon \  ({\cal H}/{\cal I})^{\otimes n+1} \otimes_{\cal H} \ad({\cal H})\stackrel{\cong}{\longrightarrow}
\big[{\cal H}^{\otimes_{\cal B} n+1}\big]_{\cal B},
\label{aux-map-psi}
\\
\psi_n\big(\big(\wbar{g^0}\otimes\cdots\otimes \wbar{g^n}\big)\otimes_{\cal H} h \big)
=\big[g^n\ps{2}hS\big(g^0\ps{1}\big) \otimes_{\cal B} g^0\ps{2}S\big(g^1\ps{1}\big) \otimes_{\cal B} \cdots \otimes_{\cal B}g^{n-1}\ps{2}S\big(g^n\ps{1}\big)\big]_{\cal B},
\nonumber
\end{gather}
with the inverse
\begin{gather}
\varphi_n\colon \ \big[{\cal H}^{\otimes_{\cal B} n+1}\big]_{{\cal B}} \stackrel{\cong}{\longrightarrow} ({\cal H}/{\cal
I})^{\otimes n+1} \otimes_{\cal H} \ad({\cal H}),
\nonumber
\\
\varphi_n\big(\big[h^0\otimes_{\cal B}\cdots\otimes_{\cal B} h^n\big]_{\cal B} \big)\nonumber\\
\qquad{}
=\big(\wbar{h^1\ps{2}\cdots h^n\ps{2}} \otimes\cdots\otimes \wbar{h^{n-1}\ps{n}h^n\ps{n}}
\otimes
\wbar{h^n\ps{n+1}} \otimes \wbar{1}\big) \otimes_{\cal H} h^0h^1\ps{1}\cdots h^n\ps{1}.
\label{aux-map-vp}
\end{gather}
\end{Lemma}

\begin{proof}
First of all, we have to prove that the maps are well def\/ined.
Let us begin with~\eqref{aux-map-vp}.
We observe that
\begin{gather*}
\varphi_n\big(\big[h^0\otimes_{\cal B}\cdots\otimes_{\cal B} h^nb\big]_{\cal B}\big)
\\
\qquad
=\big(\wbar{h^1\ps{2}\cdots (h^nb)\ps{2}}
\otimes\dots\otimes \wbar{h^{n-1}\ps{n}(h^nb)\ps{n}} \otimes
\wbar{(h^nb)\ps{n+1}} \otimes \wbar{1} \big)
\otimes_{\cal H} h^0h^1\ps{1}\cdots (h^nb)\ps{1}
\\
\qquad{}
=\big(\wbar{h^1\ps{2}\cdots h^n\ps{2}b\ps{2}}
\otimes\dots\otimes \wbar{h^{n-1}\ps{n}h^n\ps{n}b\ps{n}} \otimes
\wbar{h^n\ps{n+1}b\ps{n+1}} \otimes \wbar{1} \big)\!
\otimes_{\cal H} h^0h^1\ps{1}\cdots h^n\ps{1}b\ps{1}
\\
\qquad
=\big(\wbar{h^1\ps{2}\cdots h^n\ps{2}}b\ps{2}
\otimes\dots\otimes \wbar{h^{n-1}\ps{n}h^n\ps{n}}b\ps{n} \otimes
\wbar{h^n\ps{n+1}}b\ps{n+1} \otimes \wbar{1}\big) \otimes_{\cal H} h^0h^1\ps{1}\cdots h^n\ps{1}b\ps{1}
\\
\qquad
\stackrel{\eqref{coinviter}}{=}
\big(\wbar{h^1\ps{2}\cdots h^n\ps{2}}b\ps{2}
\otimes\dots\otimes \wbar{h^{n-1}\ps{n}h^n\ps{n}}b\ps{n} \otimes
\wbar{h^n\ps{n+1}}b\ps{n+1} \otimes \wbar{b\ps{n+2}}\big)\\
\qquad\hphantom{\stackrel{\eqref{coinviter}}{=}}{}
\otimes_{\cal H} h^0h^1\ps{1}\cdots h^n\ps{1}b\ps{1}
\\
\qquad{}
=\big(\wbar{h^1\ps{2}\cdots h^n\ps{2}}b\ps{2}
\otimes\dots\otimes \wbar{h^{n-1}\ps{n}h^n\ps{n}}b\ps{n} \otimes
\wbar{h^n\ps{n+1}}b\ps{n+1} \otimes \wbar{1}b\ps{n+2}\big)\\
\qquad{}
\hphantom{=}{}
 \otimes_{\cal H} h^0h^1\ps{1}\cdots h^n\ps{1}b\ps{1}
\\
\qquad
=\big(\wbar{h^1\ps{2}\cdots h^n\ps{2}}
\otimes\dots\otimes \wbar{h^{n-1}\ps{n}h^n\ps{n}} \otimes \wbar{h^n\ps{n+1}}
\otimes \wbar{1}\big)b\ps{2} \otimes_{\cal H} h^0h^1\ps{1}\cdots h^n\ps{1}b\ps{1}
\\
\qquad
=\big(\wbar{h^1\ps{2}\cdots h^n\ps{2}}
\otimes\dots\otimes \wbar{h^{n-1}\ps{n}h^n\ps{n}} \otimes \wbar{h^n\ps{n+1}}
\otimes \wbar{1}\big) \otimes_{\cal H} b\ps{2}\triangleright (h^0h^1\ps{1}\cdots h^n\ps{1}b\ps{1})
\\
\qquad{}
\stackrel{\eqref{HonH}}{=}
\big(\wbar{h^1\ps{2}\cdots h^n\ps{2}}
\otimes\dots\otimes \wbar{h^{n-1}\ps{n}h^n\ps{n}} \otimes \wbar{h^n\ps{n+1}}
\otimes \wbar{1}\big) \otimes_{\cal H} bh^0h^1\ps{1}\cdots h^n\ps{1}
\\
\qquad{}
=\varphi_n\big(\big[bh^0\otimes_{\cal B}\cdots\otimes_{\cal B} h^n\big]_{\cal B}\big),
\end{gather*}
that
\begin{gather*}
\varphi_n\big(\big[h^0\otimes_{\cal B} bh^{1}\otimes_{\cal B}\cdots\otimes_{\cal B} h^n\big]_{\cal B}\big)
\\
\qquad
=\big(\wbar{(bh^1)\ps{2}\cdots h^n\ps{2}} \otimes\cdots\otimes \wbar{h^{n-1}\ps{n}h^n\ps{n}} \otimes \wbar{h^n\ps{n+1}}
\otimes \wbar{1}\big) \otimes_{\cal H} h^0(bh^1)\ps{1}\cdots h^n\ps{1}
\\
\qquad
=\big(\wbar{b\ps{2}h^1\ps{2}\cdots h^n\ps{2}} \otimes\cdots\otimes \wbar{h^{n-1}\ps{n}h^n\ps{n}} \otimes
\wbar{h^n\ps{n+1}} \otimes \wbar{1}\big) \otimes_{\cal H} h^0b\ps{1}h^1\ps{1}\cdots h^n\ps{1}
\\
\qquad
=\big(\wbar{b\ps{2}}h^1\ps{2}\cdots h^n\ps{2} \otimes\cdots\otimes \wbar{h^{n-1}\ps{n}h^n\ps{n}} \otimes
\wbar{h^n\ps{n+1}} \otimes \wbar{1}\big) \otimes_{\cal H} h^0b\ps{1}h^1\ps{1}\cdots h^n\ps{1}
\\
\qquad
\stackrel{\eqref{coinv}}{=}
\big(\wbar{1}h^1\ps{2}\cdots h^n\ps{2} \otimes\cdots\otimes \wbar{h^{n-1}\ps{n}h^n\ps{n}} \otimes \wbar{h^n\ps{n+1}}
\otimes \wbar{1}\big) \otimes_{\cal H} h^0bh^1\ps{1}\cdots h^n\ps{1}
\\
\qquad
=\big(\wbar{h^1\ps{2}\cdots h^n\ps{2}} \otimes\cdots\otimes \wbar{h^{n-1}\ps{n}h^n\ps{n}} \otimes \wbar{h^n\ps{n+1}}
\otimes \wbar{1}\big) \otimes_{\cal H} h^0bh^1\ps{1}\cdots h^n\ps{1}
\\
\qquad
=\varphi_n\big(\big[h^0b\otimes_{\cal B} h^{1}\otimes_{\cal B}\cdots\otimes_{\cal B} h^n\big]_{\cal B}\big),
\end{gather*}
and that
\begin{gather*}
\varphi_n\big(\big[h^0\otimes_{\cal B}\cdots\otimes_{\cal B} h^i\otimes_{\cal B} bh^{i+1}\otimes_{\cal
B}\dots\otimes_{\cal B} h^n\big]_{\cal B}\big)
\\
\qquad
=\big(\wbar{h^1\ps{2}\cdots h^i\ps{2}(bh^{i+1})\ps{2}\cdots h^n\ps{2}}\otimes\dots\otimes
\wbar{h^i\ps{i+1}(bh^{i+1})\ps{i+1}\cdots h^n\ps{i+1}}
\\
\qquad\phantom{=}{}\otimes
\wbar{(bh^{i+1})\ps{i+2}\cdots h^n\ps{i+2}} \otimes \wbar{h^{i+2}\ps{i+3}\cdots
h^n\ps{i+3}}\otimes\cdots\otimes\wbar{h^n\ps{n+1}} \otimes \wbar{1}\big) \\
 \qquad\phantom{=}{}\otimes_{\cal H}
h^0h^1\ps{1}\cdots h^{i}\ps{1}\big(bh^{i+1}\big)\ps{1}\cdots h^n\ps{1}
\\
\qquad
=\big(\wbar{h^1\ps{2}\cdots h^i\ps{2}b\ps{2}h^{i+1}\ps{2}\cdots h^n\ps{2}}\otimes\dots\otimes
\wbar{h^i\ps{i+1}b\ps{i+1}h^{i+1}\ps{i+1}\cdots h^n\ps{i+1}}
\\
\qquad\phantom{=}{} \otimes
\wbar{b\ps{i+2}h^{i+1}\ps{i+2}\cdots h^n\ps{i+2}} \otimes \wbar{h^{i+2}\ps{i+3}\cdots
h^n\ps{i+3}}\otimes\cdots\otimes\wbar{h^n\ps{n+1}} \otimes \wbar{1}\big) \\
\qquad\phantom{=}{} \otimes_{\cal H}
h^0h^1\ps{1}\cdots h^{i}\ps{1}b\ps{1}h^{i+1}\ps{1}\cdots h^n\ps{1}
\\
\qquad
=\big(\wbar{h^1\ps{2}\cdots h^i\ps{2}b\ps{2}h^{i+1}\ps{2}\cdots h^n\ps{2}}\otimes\dots\otimes
\wbar{h^i\ps{i+1}b\ps{i+1}h^{i+1}\ps{i+1}\cdots h^n\ps{i+1}}
\\
\qquad\phantom{=}{} \otimes
\wbar{b\ps{i+2}}h^{i+1}\ps{i+2}\cdots h^n\ps{i+2} \otimes \wbar{h^{i+2}\ps{i+3}\cdots
h^n\ps{i+3}}\otimes\cdots\otimes\wbar{h^n\ps{n+1}} \otimes \wbar{1}\big) \\
\qquad\phantom{=}{}\otimes_{\cal H}
h^0h^1\ps{1}\cdots h^{i}\ps{1}b\ps{1}h^{i+1}\ps{1}\cdots h^n\ps{1}
\\
\qquad
\stackrel{\eqref{coinviter}}{=}
\big(\wbar{h^1\ps{2}\cdots h^i\ps{2}b\ps{2}h^{i+1}\ps{2}\cdots h^n\ps{2}}\otimes\dots\otimes
\wbar{h^i\ps{i+1}b\ps{i+1}h^{i+1}\ps{i+1}\cdots h^n\ps{i+1}}
\\
\qquad\phantom{=}{}
\otimes
\wbar{1}h^{i+1}\ps{i+2}\cdots h^n\ps{i+2} \otimes \wbar{h^{i+2}\ps{i+3}\cdots
h^n\ps{i+3}}\otimes\cdots\otimes\wbar{h^n\ps{n+1}} \otimes \wbar{1}\big) \\
\qquad\phantom{=}{} \otimes_{\cal H}
h^0h^1\ps{1}\cdots h^{i}\ps{1}b\ps{1}h^{i+1}\ps{1}\cdots h^n\ps{1}
\\
\qquad{}
=\big(\wbar{h^1\ps{2}\cdots (h^ib)\ps{2}h^{i+1}\ps{2}\cdots h^n\ps{2}}\otimes\cdots\otimes
\wbar{(h^ib)\ps{i+1}h^{i+1}\ps{i+1}\cdots h^n\ps{i+1}}
\\
\qquad\phantom{=}{} \otimes
\wbar{h^{i+1}\ps{i+2}\cdots h^n\ps{i+2}} \otimes \wbar{h^{i+2}\ps{i+3}\cdots
h^n\ps{i+3}}\otimes\cdots\otimes\wbar{h^n\ps{n+1}} \otimes \wbar{1}\big) \\
\qquad\phantom{=}{} \otimes_{\cal H}
h^0h^1\ps{1}\cdots (h^{i}b)\ps{1}h^{i+1}\ps{1}\cdots h^n\ps{1}
\\
\qquad{}
=\varphi_n\big(\big[h^0\otimes_{\cal B}\cdots\otimes_{\cal B} h^ib\otimes_{\cal B} h^{i+1}\otimes_{\cal
B}\cdots\otimes_{\cal B} h^n\big]_{\cal B}\big).
\end{gather*}

As for~\eqref{aux-map-psi}, we f\/irst note that by~\eqref{trans} the right hand side depends only on the representatives
appearing on the left hand side.
Moreover, for any $p\in{\cal H}$,
\begin{gather*}
\psi_n\big(\big(\wbar{g^0}\otimes\cdots\otimes \wbar{g^n}\big)\cdot p\otimes_{\cal H} h \big)
=\psi_n\big(\big(\wbar{g^0p\ps{1}}\otimes\cdots\otimes \wbar{g^np\ps{n+1}}\big)\otimes_{\cal H} h \big)
\\
\qquad
=\big[g^n\ps{2}p\ps{n+1}hS\big(g^0\ps{1}p\ps{1}\big) \otimes_{\cal B} g^0\ps{2}p\ps{2}S\big(g^1\ps{1}p\ps{3}\big) \otimes_{\cal B} \cdots\\
\qquad\phantom{=}{}
\otimes_{\cal B} g^{n-1}\ps{2}p\ps{2n}S\big(g^n\ps{1}p\ps{2n+1})\big]_{\cal B}
\\
\qquad
=\big[g^n\ps{2}p\ps{n+1}hS(p\ps{1})S\big(g^0\ps{1}\big) \otimes_{\cal B} g^0\ps{2}p\ps{2}S(p\ps{3})S\big(g^1\ps{1}\big) \otimes_{\cal B}
\cdots
\\
\qquad\phantom{=}{}
\otimes_{\cal B} g^{n-1}\ps{2}p\ps{2n}S(p\ps{2n+1})S\big(g^n\ps{1}\big)\big]_{\cal B}
\\
\qquad
=\big[g^n\ps{2}p\ps{2}hS(p\ps{1})S\big(g^0\ps{1}\big) \otimes_{\cal B} g^0\ps{2}S\big(g^1\ps{1}\big) \otimes_{\cal B} \cdots \otimes_{\cal B}
g^{n-1}\ps{2}S\big(g^n\ps{1}\big)\big]_{\cal B}
\\
\qquad
=\psi_n\big(\big(\wbar{g^0}\otimes\cdots\otimes \wbar{g^n}\big)\otimes_{\cal H} p\triangleright h \big).
\end{gather*}
Accordingly, $\psi_n$ is well-def\/ined for any $n\geq 0$.

Finally, we prove that $\varphi_n$ and $\psi_n$ are inverses to each other for any $n\geq 0$.
On one hand we have
\begin{gather*}
(\psi_n \circ \varphi_n)\big(\big[h^0\otimes_{\cal B}\cdots\otimes_{\cal B} h^n\big]_{\cal B} \big)
\\
\qquad
=\psi_n\big(\big(\wbar{h^1\ps{2}\cdots h^n\ps{2}} \otimes\cdots\otimes \wbar{h^{n-1}\ps{n}h^n\ps{n}} \otimes
\wbar{h^n\ps{n+1}} \otimes \wbar{1}\big) \otimes_{\cal H} h^0h^1\ps{1}\cdots h^n\ps{1} \big)
\\
\qquad
=\big[1\ps{2}h^0h^1\ps{1}\cdots h^n\ps{1}S\big(\big(h^1\ps{2}\cdots h^n\ps{2}\big)\ps{1}\big)\otimes_{\cal B} \big(h^1\ps{2}\cdots
h^n\ps{2}\big)\ps{2}S\big(\big(h^2\ps{3}\cdots h^n\ps{3}\big)\ps{1}\big)
\\
\qquad\phantom{=}{}
\otimes_{\cal B}\cdots \otimes_{\cal B} \big(h^{n-1}\ps{n}h^n\ps{n}\big)\ps{2}S\big(h^n\ps{n+1}\ps{1}\big)\otimes_{\cal B}
h^n\ps{n+1}\ps{2}S(1\ps{1})\big]_{\cal B}
\\
\qquad
=\big[h^0h^1\ps{1}\cdots h^n\ps{1}S\big(\big(h^1\ps{2}\cdots h^n\ps{2}\big)\ps{1}\big)\otimes_{\cal B} \big(h^1\ps{2}\cdots
h^n\ps{2}\big)\ps{2}S\big(\big(h^2\ps{3}\cdots h^n\ps{3}\big)\ps{1}\big)
\\
\qquad\phantom{=}{}
\otimes_{\cal B}\cdots \otimes_{\cal B} \big(h^{n-1}\ps{n}h^n\ps{n}\big)\ps{2}S\big(h^n\ps{n+1}\ps{1}\big)\otimes_{\cal B}
h^n\ps{n+1}\ps{2}\big]_{\cal B}
\\
\qquad
=\big[h^0h^1\ps{1}\cdots h^n\ps{1}S\big(h^1\ps{2}\ps{1}\cdots h^n\ps{2}\ps{1}\big)\otimes_{\cal B} h^1\ps{2}\ps{2}\cdots
h^n\ps{2}\ps{2}S\big(h^2\ps{3}\ps{1}\cdots h^n\ps{3}\ps{1}\big)
\\
\qquad\phantom{=}{}
\otimes_{\cal B}\cdots \otimes_{\cal B} h^{n-1}\ps{n}\ps{2}h^n\ps{n}\ps{2}S\big(h^n\ps{n+1}\ps{1}\big)\otimes_{\cal B}
h^n\ps{n+1}\ps{2}\big]_{\cal B}
\\
\qquad
=\big[h^0\cdot h^1\ps{1}\cdots h^n\ps{1}S\big(h^n\ps{2}\big)\cdots S\big(h^1\ps{2}\big)\otimes_{\cal B} h^1\ps{3}\cdot
h^2\ps{3}\cdots h^n\ps{3}S\big(h^n\ps{3}\big)\cdots S\big(h^2\ps{3}\big)
\\
\qquad\phantom{=}{}
\otimes_{\cal B}\cdots \otimes_{\cal B} h^{n-1}\ps{2n-1}h^n\ps{2n-1}S\big(h^n\ps{2n}\big)\otimes_{\cal B} h^n\ps{2n+1}\big]_{\cal B}
\\
\qquad
=\big[h^0\cdot h^1\ps{1}\cdots h^n\ps{1}S\big(h^n\ps{2}\big)\cdots S\big(h^1\ps{2}\big)\otimes_{\cal B} h^1\ps{3}\cdot
h^2\ps{3}\cdots h^n\ps{3}S\big(h^n\ps{3}\big)\cdots S\big(h^2\ps{3}\big)
\\
\qquad\phantom{=}{}
\otimes_{\cal B}\cdots \otimes_{\cal B} h^{n-1}\ps{2n-1}h^n\ps{2n-1}S\big(h^n\ps{2n}\big)\otimes_{\cal B} h^n\ps{2n+1}\big]_{\cal B}
\\
\qquad
=\big[h^0\cdot h^1\ps{1}\cdots h^{n-1}\ps{1}S\big(h^{n-1}\ps{2}\big)\cdots S\big(h^1\ps{2}\big)
\\
\qquad\phantom{=}{}
\otimes_{\cal B}
h^1\ps{3}\cdot h^2\ps{3}\cdots h^{n-1}\ps{3}S\big(h^{n-1}\ps{3}\big)\cdots S\big(h^2\ps{3}\big)\otimes_{\cal B}
\cdots \otimes_{\cal B} h^{n-1}\ps{2n-1}\otimes_{\cal B} h^n]_{\cal B} =\cdots
\\
\qquad{}
=\big[h^0\otimes_{\cal B}\cdots\otimes_{\cal B} h^n\big]_{\cal B},
\end{gather*}
while on the other hand
\begin{gather*}
(\varphi_n \circ \psi_n)\big((\wbar{g^0}\otimes\cdots\otimes \wbar{g^n})\otimes_{\cal H} h \big)
\\
\qquad
=\varphi_n\big([g^n\ps{2}hS(g^0\ps{1}) \otimes_{\cal B} g^0\ps{2}S(g^1\ps{1}) \otimes_{\cal B} \cdots \otimes_{\cal B}
g^{n-1}\ps{2}S(g^n\ps{1})]_{\cal B} \big)
\\
\qquad
=\big(\wbar{g^0\ps{2}\ps{2}S(g^1\ps{1})\ps{2}g^1\ps{2}\ps{2}S(g^2\ps{1})\ps{2} \cdots
g^{n-1}\ps{2}\ps{2}S(g^n\ps{1})\ps{2}} \otimes \cdots
\\
\qquad
\phantom{=}{}
\otimes \wbar{g^{n-2}\ps{2}\ps{n}S(g^{n-1}\ps{1})\ps{n}g^{n-1}\ps{2}\ps{n}S(g^n\ps{1})\ps{n}}\otimes
\wbar{g^{n-1}\ps{2}\ps{n+1}S(g^n\ps{1})\ps{n+1}} \otimes \wbar{1}\big)
\\
\qquad\phantom{=}{}\otimes_{\cal H}
g^n\ps{2}hS(g^0\ps{1})g^0\ps{2}\ps{1}S(g^1\ps{1})\ps{1}g^1\ps{2}\ps{1}S(g^2\ps{1})\ps{1}\cdots
g^{n-1}\ps{2}\ps{1}S(g^n\ps{1})\ps{1}
\\
\qquad
=\big(\wbar{g^0\ps{3}S(g^1\ps{1})g^1\ps{4}S(g^2\ps{2}) \cdots g^{n-1}\ps{n+2}S(g^n\ps{n})} \otimes \cdots
\\
\qquad
\phantom{=}
\otimes \wbar{g^{n-2}\ps{2n-1}S(g^{n-1}\ps{1})g^{n-1}\ps{2n}S(g^n\ps{2})}\otimes
\wbar{g^{n-1}\ps{2n+1}S(g^n\ps{1})} \otimes \wbar{1}\big)
\\
\qquad\phantom{=}\otimes_{\cal H}
g^n\ps{n+2}hS(g^0\ps{1})g^0\ps{2}S(g^1\ps{2})g^1\ps{3}S(g^2\ps{3})\cdots g^{n-1}\ps{n+1}S(g^n\ps{n+1})
\\
\qquad
=\big(\wbar{g^0S(g^n\ps{n})} \otimes \cdots\otimes \wbar{g^{n-2}S(g^n\ps{2})}\otimes \wbar{g^{n-1}S(g^n\ps{1})} \otimes
\wbar{1}\big)\otimes_{\cal H} g^n\ps{n+2}hS(g^n\ps{n+1})
\\
\qquad
=\big(\wbar{g^0S(g^n\ps{1}\ps{n})} \otimes\cdots\otimes \wbar{g^{n-1}S(g^n\ps{1}\ps{1})}\otimes
\wbar{1}\big)\otimes_{\cal H} g^n\ps{2}hS(g^n\ps{1}\ps{n+1})
\\
\qquad
=\big(\big(\wbar{g^0} \otimes\cdots\otimes \wbar{g^{n-1}}\big)\cdot S\big(g^n\ps{1}\big)\otimes \wbar{1}\big)\otimes_{\cal H}
g^n\ps{2} \triangleright h
\\
\qquad
=\big(\big(\wbar{g^0} \otimes\cdots\otimes \wbar{g^{n-1}}\big)\cdot S(g^n\ps{1})\otimes \wbar{1}\big)\cdot g^n\ps{2}
\otimes_{\cal H} h
\\
\qquad
=\big(\wbar{g^0} \otimes\cdots\otimes \wbar{g^{n-1}}\otimes \wbar{g^n}\big) \otimes_{\cal H} h. \tag*{\qed}
\end{gather*}
\renewcommand{\qed}{}
\end{proof}

We are now ready to state our main result.
\begin{Theorem}
\label{dirpict}
Let ${\cal I}\subseteq {\cal H}$ be a~coideal right ideal in a~Hopf algebra ${\cal H}$
such that ${\cal H}^{\co{\cal H}/{\cal I}}\subseteq {\cal H}$ is an ${\cal H}/{\cal I}$-Galois extension.
Let also $\ad({\cal H})={\cal H}$ be the left-right SAYD module with the right adjoint action, and the left
coaction given by the comultiplication of~${\cal H}$.
Then there exists an isomorphism
\begin{gather*}
\psi_{n}\colon \  {\rm C}_{n}({\cal H}/{\cal I},\ad({\cal H}))_{{\cal H}}\longrightarrow {\rm C}_{n}\big({\cal H}\mid{\cal H}^{\co{\cal H}/{\cal I}}\big)
\end{gather*}
of cyclic modules, defined by~\eqref{aux-map-psi}, \eqref{aux-map-vp}.
\end{Theorem}
\begin{proof}
Let us, as above, adopt the notation ${\cal B}:={\cal H}^{\co{\cal H}/{\cal I}}$.
We shall f\/irst go through the commutation with the face operators.
For $i=0$,
\begin{gather*}
\varphi_{n-1}d_0\big(\big[h^0\otimes_{\cal B}\cdots\otimes_{\cal B} h^n\big]_{\cal B}\big) =\varphi_{n-1}\big(\big[h^0h^1\otimes_{\cal
B}\cdots\otimes_{\cal B} h^n\big]_{\cal B}\big)
\\
\qquad
=\big(\wbar{h^2\ps{2}\cdots h^n\ps{2}} \otimes\cdots\otimes \wbar{h^{n-1}\ps{n-1}h^n\ps{n-1}} \otimes \wbar{h^n\ps{n}}
\otimes \wbar{1}\big) \otimes_{\cal H} h^0h^1h^2\ps{1}\cdots h^n\ps{1}
\\
\qquad
=d_0\varphi_{n}\big(\big[h^0\otimes_{\cal B}\cdots\otimes_{\cal B} h^n\big]_{\cal B}\big).
\end{gather*}
For $1\leq i \leq n-1$, we observe that
\begin{gather*}
\varphi_{n-1}d_i\big(\big[h^0\otimes_{\cal B}\cdots\otimes_{\cal B} h^n\big]_{\cal B}\big) =\varphi_{n-1}([h^0\otimes_{\cal
B}\cdots\otimes_{\cal B} h^ih^{i+1}\otimes_{\cal B}\cdots\otimes_{\cal B} h^n]_{\cal B})
\\
\qquad
=\big(\wbar{h^1\ps{2}\cdots(h^ih^{i+1})\ps{2}\cdots h^n\ps{2}} \otimes\cdots\otimes \wbar{(h^ih^{i+1})\ps{i+1}\dots
h^n\ps{i+1}}
\\
\qquad
\phantom{=}{}\otimes\wbar{h^{i+2}\ps{i+2}\cdots h^n\ps{i+2}}\otimes\cdots
\otimes\wbar{h^{n-1}\ps{n-1}h^n\ps{n-1}} \otimes \wbar{h^n\ps{n}} \otimes \wbar{1}\big)\\
\qquad
\phantom{=}{}
 \otimes_{\cal H}
h^0h^1\ps{1}\cdots(h^ih^{i+1})\ps{1}\cdots h^n\ps{1}
\\
\qquad
=\big(\wbar{h^1\ps{2}\cdots h^n\ps{2}} \otimes\cdots\otimes \wbar{h^i\ps{i+1}h^{i+1}\ps{i+1}\cdots h^n\ps{i+1}}\otimes
\varepsilon\big(\wbar{h^{i+1}\ps{i+2}\cdots h^n\ps{i+2}}\big) \otimes\cdots
\\
\qquad
\phantom{=}{}
\otimes \wbar{h^{n-1}\ps{n}h^n\ps{n}} \otimes \wbar{h^n\ps{n+1}} \otimes \wbar{1}\big) \otimes_{\cal H}
h^0h^1\ps{1}\cdots h^n\ps{1}
\\
\qquad
=d_i\varphi_{n}\big(\big[h^0\otimes_{\cal B}\cdots\otimes_{\cal B} h^n\big]_{\cal B}\big).
\end{gather*}
Finally for the last face operator we have
\begin{gather*}
\varphi_{n-1}d_n\big(\big[h^0\otimes_{\cal B}\cdots\otimes_{\cal B} h^n\big]_{\cal B}\big) =\varphi_{n-1}\big(\big[h^nh^0\otimes_{\cal
B}\cdots\otimes_{\cal B} h^{n-1}\big]_{\cal B}\big)
\\
\qquad
=\big(\wbar{h^1\ps{2}\cdots h^{n-1}\ps{2}} \otimes \cdots
\otimes \wbar{h^{n-2}\ps{n-1}h^{n-1}\ps{n-1}} \otimes \wbar{h^{n-1}\ps{n}} \otimes
\wbar{1}\big)
 \otimes_{\cal H} h^nh^0h^1\ps{1}\cdots h^{n-1}\ps{1}
\\
\qquad
\stackrel{\eqref{HonH}}{=}
\big(\wbar{h^1\ps{2}\cdots h^{n-1}\ps{2}} \otimes\cdots
\otimes \wbar{h^{n-2}\ps{n-1}h^{n-1}\ps{n-1}} \otimes \wbar{h^{n-1}\ps{n}} \otimes
\wbar{1}\big) \\
\qquad\hphantom{=}{}
\otimes_{\cal H} h^n\ps{2}\triangleright \big(h^0h^1\ps{1}\cdots h^{n-1}\ps{1}h^n\ps{1}\big)
\\
\qquad
=\big(\wbar{h^1\ps{2}\cdots h^{n-1}\ps{2}} \otimes\cdots
\otimes \wbar{h^{n-2}\ps{n-1}h^{n-1}\ps{n-1}} \otimes \wbar{h^{n-1}\ps{n}} \otimes \wbar{1}\big)\cdot
h^n\ps{2} \\
\qquad
\phantom{=}
\otimes_{\cal H} h^0h^1\ps{1}\cdots h^{n-1}\ps{1}h^n\ps{1}
\\
\qquad
=\big(\wbar{h^1\ps{2}\cdots h^{n}\ps{2}} \otimes\dots\otimes \wbar{h^{n-1}\ps{n}h^{n}\ps{n}} \otimes
\wbar{h^{n}\ps{n+1}} \big) \otimes_{\cal H} h^0h^1\ps{1}\cdots h^{n}\ps{1}
\\
\qquad
=\big(\wbar{h^1\ps{2}\cdots h^{n}\ps{2}} \otimes\cdots\otimes \wbar{h^{n-1}\ps{n}h^{n}\ps{n}} \otimes
\wbar{h^{n}\ps{n+1}}\otimes\varepsilon\big(\wbar{1}\big) \big) \otimes_{\cal H} h^0h^1\ps{1}\cdots h^{n}\ps{1}
\\
\qquad
=d_n\varphi_{n}\big(\big[h^0\otimes_{\cal B}\cdots\otimes_{\cal B} h^n\big]_{\cal B}\big).
\end{gather*}
We next investigate the interaction with the degeneracy operators.
To this end, for $0\leq j \leq n-1$ we observe that
\begin{gather*}
\varphi_{n+1}s_j\big(\big[h^0\otimes_{\cal B}\cdots\otimes_{\cal B} h^n\big]_{\cal B}\big)
=\varphi_{n+1}\big(\big[h^0\otimes_{\cal B}\cdots\otimes_{\cal B} h^j\otimes_{\cal B} 1\otimes_{\cal B} h^{j+1}\otimes_{\cal
B}\cdots\otimes_{\cal B} h^{n}\big]_{\cal B}\big)
\\
\qquad
=\big(\wbar{h^1\ps{2}\cdots h^n\ps{2}} \otimes\cdots
\otimes \wbar{h^j\ps{j+1}\cdots h^n\ps{j+1}} \otimes \wbar{h^{j+1}\ps{j+2}\cdots h^n\ps{j+2}}\\
\qquad
\phantom{=}{}
\otimes
\wbar{h^{j+1}\ps{j+3}\dots h^n\ps{j+3}}\otimes\cdots
\\
\qquad
\phantom{=}{}
\otimes \wbar{h^{n-1}\ps{n+1}h^n\ps{n+1}} \otimes \wbar{h^n\ps{n+2}} \otimes \wbar{1}\big) \otimes_{\cal H}
h^0h^1\ps{1}\cdots h^n\ps{1}
\\
\qquad
=\big(\wbar{h^1\ps{2}\cdots h^n\ps{2}} \otimes\cdots\otimes \wbar{h^j\ps{j+1}\cdots h^n\ps{j+1}} \otimes
\Delta\big(\wbar{h^{j+1}\ps{j+2}\cdots h^n\ps{j+2}}\big)\otimes \cdots
\\
\qquad
\phantom{=}{}
\otimes \wbar{h^{n-1}\ps{n}h^n\ps{n}} \otimes \wbar{h^n\ps{n+1}} \otimes \wbar{1}\big) \otimes_{\cal H}
h^0h^1\ps{1}\cdots h^n\ps{1}
\\
\qquad
=s_j\varphi_{n}\big(\big[h^0\otimes_{\cal B}\cdots\otimes_{\cal B} h^n\big]_{\cal B}\big).
\end{gather*}
Let us f\/inally check the cyclic operators.
To this end we have
\begin{gather*}
\varphi_{n}t_n\big(\big[h^0\otimes_{\cal B}\cdots\otimes_{\cal B} h^n\big]_{\cal B}\big) =\varphi_{n+1}\big(\big[h^n\otimes_{\cal B}
h^0\otimes_{\cal B}\cdots\otimes_{\cal B} h^{n-1}\big]_{\cal B}\big)
\\
\qquad
=\big(\wbar{h^0\ps{2}\cdots h^{n-1}\ps{2}} \otimes\cdots\otimes \wbar{h^{n-2}\ps{n}h^{n-1}\ps{n}} \otimes
\wbar{h^{n-1}\ps{n+1}} \otimes \wbar{1}\big) \otimes_{\cal H} h^nh^0\ps{1}\cdots h^{n-1}\ps{1}
\\
\qquad
=\big(\wbar{h^0\ps{2}h^1\ps{2}\cdots h^n\ps{2}}\otimes \wbar{h^1\ps{3}\cdots h^n\ps{3}} \otimes\cdots
\\
\qquad
\phantom{=}{}
\otimes \wbar{h^{n-1}\ps{n+1}h^n\ps{n+1}} \otimes \wbar{h^n\ps{n+2}} \big) \otimes_{\cal H}
h^0\ps{1}h^1\ps{1}\cdots h^n\ps{1}
\\
\qquad
=\big(\wbar{1}\cdot (h^0h^1\ps{1}\cdots h^n\ps{1})\ps{2}\otimes \wbar{h^1\ps{2}\cdots h^n\ps{2}} \otimes\cdots
\\
\qquad
\phantom{=}{}
\otimes \wbar{h^{n-1}\ps{n}h^n\ps{n}} \otimes \wbar{h^n\ps{n+1}} \big) \otimes_{\cal H}
(h^0h^1\ps{1}\cdots h^n\ps{1})\ps{1}
\\
\qquad
=t_n\varphi_{n}\big(\big[h^0\otimes_{\cal B}\cdots\otimes_{\cal B} h^n\big]_{\cal B}\big). \tag*{\qed}
\end{gather*}
\renewcommand{\qed}{}
\end{proof}

\subsection{Comparison with the Jara--\c{S}tefan isomorphism\\ in the homogeneous quotient-Hopf--Galois case}

In this subsection we compare the isomorphism of Theorem~\ref{dirpict} with that of~\cite[Theorem~3.7]{JaraStef06} in
the case of the homogeneous quotient Hopf--Galois extensions.

We note that in the case of ${\cal I}$ being a~Hopf ideal of ${\cal H}$, analogously to~\eqref{aux-map-vp} we have a~map
\begin{gather}
\label{jara-stefan}
\overline{\varphi}_n\colon \ \big[{\cal H}^{\otimes_{\cal B} n+1}\big]_{{\cal B}} \stackrel{\cong}{\longrightarrow} ({\cal
H}/{\cal I})^{\otimes n+1} \otimes_{{\cal H}/{\cal I}} [\ad({\cal H})]_{{\cal B}}
\\
\overline{\varphi}_n\big([h^0\otimes_{\cal B}\cdots\otimes_{\cal B} h^n]_{\cal B} \big)\nonumber
\\
\qquad
=\big(\wbar{h^1\ps{2}\cdots h^n\ps{2}} \otimes\cdots\otimes \wbar{h^{n-1}\ps{n}h^n\ps{n}} \otimes
\wbar{h^n\ps{n+1}} \otimes \wbar{1}\big) \otimes_{{\cal H}/{\cal I}} \big[h^0h^1\ps{1}\cdots h^n\ps{1}\big]_{{\cal B}},\nonumber
\end{gather}
which is well def\/ined by the Hopf ideal assumption.
It is easy to see that it is a~special case of the isomorphism introduced by Jara--\c{S}tefan~\cite[Theorem~3.7]{JaraStef06}
originally in the context of Hopf--Galois extensions, when restricted to the context of homogeneous
quotient-Hopf--Galois extensions.

Compared to~\eqref{aux-map-vp}, on the right hand side we tensorize over ${\cal H}$ from the right by $\ad({\cal H})$
(with the left ${\cal H}$-module structure~\eqref{HonM}) while in~\eqref{jara-stefan}
the tensor product by $[\ad({\cal H})]_{{\cal B}}$ is over ${\cal H}/{\cal I}$ (with the left Miyashita--Ulbrich ${\cal H}/{\cal I}$-module structure).

However, in the case of homogeneous quotient-Hopf--Galois extensions the isomorphisms~\eqref{aux-map-vp}
and~\eqref{jara-stefan} coincide.
Indeed, in this case ${\cal H}/{\cal I}$ is generated by $\overline{1}$ both as a~left and a~right ${\cal H}$-module,
with respect to the multiplication from the left or from the right in ${\cal H}$.
Accordingly, the degree zero component of the isomorphism~\eqref{aux-map-vp} yields
\begin{gather*}
\varphi_{0}\colon \   [{\cal H} ]_{{\cal B}} \stackrel{\cong}{\longrightarrow} ({\cal H}/{\cal I})\otimes_{{\cal H}}
\ad({\cal H}),\qquad \varphi_0\big([h]_{\cal B} \big) =\wbar{1} \otimes_{{\cal H}} h,
\end{gather*}
which is left ${\cal H}/{\cal I}$-linear with respect to the Miyashita--Ulbrich action on $[{\cal H}]_{{\cal
B}}$ on the left hand side, and the algebra map ${\cal H}\rightarrow {\cal H}/{\cal I}$ on the right hand side, since
\begin{gather*}
\varphi_0\big(\wbar{h}\triangleright [h']_{\cal B} \big)
=\varphi_0\big([h\triangleright h']_{\cal B} \big)
=\wbar{1}\otimes_{{\cal H}} h\triangleright h'
=\wbar{1}\cdot h \otimes_{{\cal H}} \triangleright h'
\\
\phantom{\varphi_0\big(\wbar{h}\triangleright [h']_{\cal B} \big)}
=\wbar{h} \otimes_{{\cal H}} h' =h\cdot\wbar{1} \otimes_{{\cal H}} h'=h\cdot\varphi_0\big([h']_{\cal B} \big).
\end{gather*}
We thus obtain an isomorphism of functors
\begin{gather*}
-\otimes_{{\cal H}/{\cal I}}\left[\ad({\cal H})\right]_{{\cal B}}\cong -\otimes_{{\cal H}/{\cal I}}\big(({\cal
H}/{\cal I})\otimes_{{\cal H}}\ad({\cal H})\big)=-\otimes_{{\cal H}}\ad({\cal H}).
\end{gather*}

\subsection{A spectral sequence}

We note that the right hand side of Theorem~\ref{dirpict}, resp.\
the left hand side of Theorem~\ref{dualpict}, compute the homology of the algebra extension ${\cal H}$ relative to
${\cal B}$, resp.\
coextension ${\cal H}$ corelative to ${\cal C}$, while the other side computes the cyclic dual homology of the
Pontryagin dual objects.
In order to be able to investigate the latter homologies, in this subsection we shall develop computational tools.

We will focus on the Hochschild homology groups of the relative homology of the extension.
\begin{Theorem}
Let ${\cal I}\subseteq {\cal H}$ be a~coideal right ideal in a~Hopf algebra ${\cal H}$
such that ${\cal H}^{\co{\cal H}/{\cal I}}\subseteq {\cal H}$ be a~homogeneous ${\cal H}/{\cal I}$-Galois extension.
Then there exists a~spectral sequence $($constructed in the proof$)$ such that
\begin{gather*}
{\rm HH}_\bullet\big({\cal H}\mid {\cal H}^{\co{\cal H}/{\cal I}}\big)={\rm E}^2_{\bullet,0},
\qquad
{\rm E}^2_{\bullet,\bullet}\Longrightarrow \Tor^{{\cal H}}_{\bullet}(k, \ad({\cal H})).
\end{gather*}
In particular, we have a~five-term exact sequence
\begin{gather*}
\phantom{\cdots\to }
\Tor^{{\cal H}}_2(k,\ad({\cal H}))\to {\rm HH}_2\big({\cal H}\mid {\cal H}^{\co{\cal H}/{\cal I}}\big)\to
{\rm H}_0\big(\Tor^{{\cal H}}_1(({\cal H}/{\cal I})^{\otimes\bullet +1},\ad({\cal H}))\big) \to \cdots
\\
\cdots\to \Tor^{{\cal H}}_1(k, \ad({\cal H})) \to {\rm HH}_1\big({\cal A}\mid {\cal H}^{\co{\cal H}/{\cal I}}\big) \to 0.
\end{gather*}
\end{Theorem}

\begin{proof}
Let us consider the cyclic dual ${\cal C}_{\bullet}$ to the standard cocyclic object of the coalgebra ${\cal C}:={\cal
H}/{\cal I}$, consisting of the tensor powers of ${\cal C}$.
More precisely, ${\cal C}_{p}={\cal C}^{\otimes p+1}$, with the boundary map $\partial\colon {\cal C}_{p}\rightarrow {\cal
C}_{p-1}$ being
\begin{gather*}
\partial\big(c^0\otimes\cdots\otimes c^p\big)=\sum\limits_{i=0}^{p}(-1)^i
c^0\otimes\cdots\otimes\varepsilon\big(c^i\big)c^{i+1}\otimes\cdots\otimes c^p,
\end{gather*}
where the cyclic order of length $p+1$ is assumed.
Note that $\partial$ is a~morphism of right ${\cal H}$-modules.
The operator $h\colon {\cal C}_{p}\rightarrow {\cal C}_{p+1}$
\begin{gather*}
h\big(c^0\otimes\cdots\otimes c^p\big):=\overline{1}\otimes c^0\otimes\cdots\otimes c^p
\end{gather*}
is a~homotopy contracting this complex to~$k$ concentrated at zero degree, see for instance~\cite{KhalRang05}.
Let ${\cal M}_{\bullet}$ be a~f\/lat resolution of the left ${\cal H}$-module ${\cal M}:=\ad({\cal H})$, and consider
the total complex ${\cal C}_{\bullet}\otimes_{{\cal H}}{\cal M}_{\bullet}$.
We have two spectral sequences abutting to the total homology.
The f\/irst page of the f\/irst one reads as
\begin{gather*}
{\rm E}^1_{p,q}={\rm H}_{q}({\cal C}_{p}\otimes_{{\cal H}}{\cal M}_{\bullet})=\Tor^{{\cal H}}_q({\cal C}_{p},{\cal M}),
\end{gather*}
hence its second page is computed as
\begin{gather*}
{\rm E}^2_{p,q}={\rm H}_p\big(\Tor^{{\cal H}}_q({\cal C}_{\bullet},{\cal M})\big).
\end{gather*}
By Theorem~\ref{dirpict}, we are interested in
\begin{gather}
\label{rel}
{\rm HH}_n\big({\cal H}\mid {\cal H}^{\co{\cal H}/{\cal I}}\big)={\rm E}^2_{n,0}.
\end{gather}

Now, let us compute the f\/irst page of the second spectral sequence, the transposed analogue of the f\/irst one, abutting
to the total homology, which by f\/latness of the resolution ${\cal M}_{\bullet}$ and the acyclicity of ${\cal
C}_{\bullet}$ can be rewritten as
\begin{gather*}
^{\top }{\rm E}^1_{p,q}={\rm H}_{p}({\cal C}_{\bullet}\otimes_{{\cal H}}{\cal M}_{q})={\rm H}_{p}({\cal
C}_{\bullet})\otimes_{{\cal H}}{\cal M}_{q}=\begin{cases}
k\otimes_{{\cal H}}{\cal M}_{q}& \text{for}\quad  p=0,
\\
0 & \text{for}\quad  p>0.
\end{cases}
\end{gather*}
As a~result, degenerating at the second page,
\begin{gather*}
^{\top }{\rm E}^2_{p,q}=\begin{cases}
\Tor^{{\cal H}}_{q}(k,{\cal M})& \text{for}\quad  p=0,
\\
0 & \text{for}\quad  p>0,
\end{cases}
\end{gather*}
this spectral sequence yields the total cohomology
\begin{gather}
\label{tot}
{\rm H}_{n}({\cal C}_{\bullet}\otimes_{{\cal H}}{\cal M}_{\bullet})=\Tor^{{\cal H}}_{n}(k,{\cal M}).
\end{gather}
Finally, we use~\eqref{rel},~\eqref{tot} and the canonical homological f\/ive-term exact sequence
\begin{gather*}
{\rm H}_2\rightarrow {\rm E}^2_{2,0}\stackrel{d}{\rightarrow} {\rm E}^2_{0,1}\rightarrow {\rm H}_1\rightarrow {\rm
E}^2_{1,0}\rightarrow 0
\end{gather*}
to f\/inish the proof.
The second arrow is the boundary map of the second page of the f\/irst spectral sequence, and the next arrow is induced~by
the augmentation ${\cal C}_{\bullet}\rightarrow k$.
\end{proof}

Next we show that the above theorem generalizes the classical result from the case of the homogeneous ${\cal H}$-Galois
extension $k\subseteq{\cal H}$ to arbitrary homogeneous quotient coalgebra-Galois extensions ${\cal B}\subseteq {\cal
H}$.
The following can be regarded as an independent proof of this classical result.
\begin{Corollary}
For any Hopf algebra ${\cal H}$,
\begin{gather}
\label{aux-Hoch-Tor}
{\rm HH}_\bullet({\cal H})=\Tor^{{\cal H}}_{\bullet}(k, \ad({\cal H})),
\end{gather}
where on the right hand side the left ${\cal H}$-module structure on ${\cal H}$ comes from its canonical SAYD module
structure.
\end{Corollary}
\begin{proof}
Since for every right ${\cal H}$-module ${\cal N}$ the invertible linear map
\begin{gather*}
{\cal N}\otimes {\cal H}\longrightarrow {\cal N}\otimes {\cal H},\qquad n\otimes h\mapsto nS(h\ps{1})\otimes h\ps{2}
\end{gather*}
makes the diagonal right ${\cal H}$-module ${\cal N}\otimes {\cal H}$ free, hence f\/lat, then by induction the diagonal
right ${\cal H}$-module ${\cal C}_{\bullet}={\cal H}^{\otimes\bullet+1}$ is f\/lat.
This implies that
\begin{gather}
\label{degen}
{\rm E}^2_{p,q}={\rm H}_p\big(\Tor^{{\cal H}}_q({\cal C}_{\bullet},\ad({\cal H}))\big)=\begin{cases}
{\rm H}_p({\cal C}_{\bullet}\otimes_{{\cal H}}\ad({\cal H}))& \text{for}\quad  q=0,
\\
0 & \text{for}\quad  q>0,
\end{cases}
\end{gather}
hence the f\/irst spectral sequence degenerates at the second page as well, and therefore we obtain
\begin{gather*}
{\rm HH}_n({\cal H})={\rm HH}_n({\cal H}\mid k) \stackrel{\eqref{rel}}{=} {\rm E}^2_{n,0}\stackrel{\eqref{degen}}{=}{\rm
H}_{n}({\cal C}_{\bullet}\otimes_{{\cal H}}{\cal M}_{\bullet})\stackrel{\eqref{tot}}{=}\Tor^{{\cal H}}_{n}(k,{\cal H}).  \tag*{\qed}
\end{gather*}
\renewcommand{\qed}{}
\end{proof}

\subsection{The isomorphism of Pontryagin duals}

In this subsection, we discuss the Pontryagin dual of the results of the previous subsection (we will call the results
of Subsection~\ref{subsect-iso-cyclic-obj} \emph{direct picture}, and their dual counterparts \emph{dual picture}) in
detail.
Hence, the result of the present subsection can be interpreted as the \emph{cyclic-homological dual Takeuchi--Galois
transform} accompanying the Takeuchi--Galois transform ${\cal B}\mapsto{\cal H}/{\cal B}^{+}{\cal H}$.
More explicitly, we shall obtain an isomorphism between the Hopf-cyclic homology of the ${\cal H}$-comodule algebra
${\cal B}$ and the relative homology of a~${\cal B}$-Galois coextension ${\cal D}:={\cal H} \to {\cal H}/{\cal
B}^{+}{\cal H}=:{\cal C}$ mentioned in the following theorem.
\begin{Theorem}
\label{dualpict}
Let ${\cal B}\subseteq {\cal H}$ be a~left comodule subalgebra in a~Hopf algebra ${\cal H}$ such that ${\cal
H}\rightarrow {\cal H}/{\cal B}^{+}{\cal H}$ is a~${\cal B}$-Galois coextension, and ${\rm coad}({\cal H})={\cal H}$ be
the right-left SAYD module with the right action given by the multiplication of ${\cal H}$, and the left coadjoint
coaction.
Then there exists an isomorphism
\begin{gather*}
\gamma_{n}\colon \ {\rm C}_{n}({\cal B}, {\rm coad}({\cal H}))^{{\cal H}} \longrightarrow {\rm C}_{n}({\cal H}\mid {\cal
H}/{\cal B}^{+}{\cal H})
\end{gather*}
of cyclic modules, defined by~\eqref{aux-vp-psi}, \eqref{aux-vprime-psiprime}.
\end{Theorem}
\begin{proof}
We note that the proof of Theorem~\ref{dirpict} uses only the structure maps and the relations equivalent to the
commutativity of appropriate diagrams, therefore, is in fact diagrammatical.
Then applying the formal Pontryagin duality by reversing the arrows, interchanging the \emph{left} and \emph{right},
next applying the cyclic duality and f\/inally inverting the resulting isomorphism yields a~diagrammatical proof of the
claim.
More precisely, the isomorphism is given~by
\begin{gather}
\gamma_n\colon \ {\rm coad}({\cal H}) \Box_{\cal H} {\cal B}^{\otimes n+1}
\longrightarrow (\underset{n+1~\text{many}}{\underbrace{{\cal H} \Box_{\cal C} {\cal H} \Box_{\cal C} \cdots \Box_{\cal C} {\cal H}}})^{\cal C},
\nonumber
\\
h\otimes b^0\otimes\cdots\otimes b^n \mapsto
b^1\ps{2}\cdots b^n\ps{2}h\ps{2} \otimes\cdots\otimes b^n\ps{n+1}h\ps{n+1}\otimes b^0b^1\ps{1}\cdots b^n\ps{1}h\ps{1},
\label{aux-vp-psi}
\end{gather}
and
\begin{gather}
\gamma^{-1}_n\colon \ (\underset{n+1~\text{many}}{\underbrace{{\cal H} \Box_{\cal C} {\cal H} \Box_{\cal C} \cdots \Box_{\cal C}
{\cal H}}})^{\cal C} \longrightarrow {\rm coad}({\cal H}) \Box_{\cal H} {\cal B}^{\otimes n+1},
\nonumber
\\
h^0\otimes\cdots\otimes h^n \mapsto h^n\ps{2}\otimes h^n\ps{3}S(h^0\ps{1})\otimes h^0\ps{2}S(h^1\ps{1})
\otimes\cdots\otimes h^{n-1}\ps{2}S(h^n\ps{1}).
\label{aux-vprime-psiprime}
\end{gather}

We next show that the cyclic structure on $C_\bullet({\cal C}\mid {\cal H})$ corresponds, via the
isomorphisms~\eqref{aux-vp-psi} and~\eqref{aux-vprime-psiprime}, to the Hopf-cyclic structure on $C_\bullet({\cal
B},{\rm coad}({\cal H}))^{\cal H}$.

We f\/irst note for $0\leq i \leq n-1$ that
\begin{gather*}
\gamma_{n-1}\circ d_i\big(h\otimes b^0\otimes\cdots\otimes b^n\big) =\gamma_{n-1}\big(h\otimes b^0\otimes\cdots\otimes b^ib^{i+1}\otimes\cdots\otimes b^n\big)
\\
\qquad
=b^1\ps{2}\cdots b^n\ps{2}h\ps{2} \otimes\cdots\otimes b^i\ps{i+1}b^{i+1}\ps{i+1}b^{i+2}\ps{i+1}\cdots
b^n\ps{i+1}h\ps{i+1}
\\
\qquad\phantom{=}{}
\otimes b^{i+2}\ps{i+2}\cdots b^n\ps{i+2}h\ps{i+2}
\otimes\cdots
\otimes b^n\ps{n}h\ps{n}\otimes b^0b^1\ps{1}\cdots b^n\ps{1}h\ps{1}
\\
\qquad
=b^1\ps{2}\cdots b^n\ps{2}h\ps{2} \otimes\cdots\otimes b^i\ps{i+1}b^{i+1}\ps{i+1}b^{i+2}\ps{i+1}\cdots b^n\ps{i+1}h\ps{i+1}
\\
\qquad\phantom{=}{}\otimes
\varepsilon\big(b^{i+1}\ps{i+2}b^{i+2}\ps{i+2}\cdots b^n\ps{i+2}h\ps{i+2}\big)\otimes b^{i+2}\ps{i+3}\cdots b^n\ps{i+3}h\ps{i+3}
\otimes\cdots  \\
\qquad\phantom{=}{}
\otimes b^n\ps{n+1}h\ps{n+1}\otimes b^0b^1\ps{1}\cdots b^n\ps{1}h\ps{1}
\\
\qquad
=d_i\circ \gamma_n \big(h\otimes b^0\otimes\cdots\otimes b^n\big).
\end{gather*}
As for $i=n$, we have
\begin{gather*}
\gamma_{n-1}\circ d_n\big(h\otimes b^0\otimes\cdots\otimes b^n\big) =\gamma_{n-1}\big(b^n\ps{1}h\otimes
b^n\ps{2}b^0\otimes\cdots\otimes b^{n-1}\big)
\\
\qquad
=b^1\ps{2}\cdots b^{n-1}\ps{2}b^n\ps{1}\ps{2}h\ps{2} \otimes\cdots\otimes b^{n-1}\ps{n}b^n\ps{1}\ps{n}h\ps{n}\otimes
b^n\ps{2}b^0b^1\ps{1}\cdots b^n\ps{1}h\ps{1}
\\
\qquad
=b^1\ps{2}\cdots b^n\ps{2}h\ps{2} \otimes\cdots\otimes b^{n-1}\ps{n}b^n\ps{n}h\ps{n}\otimes b^n\ps{2}b^0b^1\ps{1}\cdots b^n\ps{1}h\ps{1}
\\
\qquad
=b^1\ps{1}\cdots b^n\ps{1}h\ps{2}\ps{2} \otimes\cdots\otimes b^{n-1}\ps{n-1}b^n\ps{n-1}h\ps{2}\ps{n}\otimes
b^n\ps{n}h\ps{3}S(h\ps{1})h\ps{2}\ps{1}\varepsilon\big(b^0\big)
\\
\qquad
=b^1\ps{1}\cdots b^n\ps{1}h\ps{1} \otimes\cdots\otimes b^{n-1}\ps{n-1}b^n\ps{n-1}h\ps{n-1}\otimes b^n\ps{n}h\ps{n}\varepsilon\big(b^0\big)
\\
\qquad
=b^1\ps{2}\cdots b^n\ps{2}h\ps{2} \otimes\cdots\otimes b^n\ps{n+1}h\ps{n+1}\otimes \varepsilon\big(b^0b^1\ps{1}\cdots b^n\ps{1}h\ps{1}\big)
\\
\qquad
=d_n\circ \gamma_n \big(h\otimes b^0\otimes\dots\otimes b^n\big),
\end{gather*}
where on the third equality we used
\begin{gather*}
h\otimes b^0\ps{1}\cdots b^n\ps{1}\otimes \varepsilon\big(b^0\ps{2}\big)\otimes\cdots\otimes b^n\ps{2} =h\ps{2}\otimes
h\ps{3}S(h\ps{1}) \otimes \varepsilon\big(b^0\big)\otimes\cdots\otimes b^n.
\end{gather*}
We next note that
\begin{gather*}
\gamma_{n+1}\circ s_j\big(h\otimes b^0\otimes\cdots\otimes b^n\big) =\gamma_{n+1}\big(h\otimes \otimes b^0\otimes\cdots\otimes
b^j\otimes 1\otimes\cdots\otimes b^n\big)
\\
\qquad
=b^1\ps{2}\cdots b^n\ps{2}h\ps{2} \otimes\cdots\otimes b^{j+1}\ps{j+1}\cdots b^n\ps{j+1}h\ps{j+1}\otimes
b^{j+1}\ps{j+2}\cdots b^n\ps{j+2}h\ps{j+2} \otimes\cdots
\\
\qquad
\phantom{=}{}
\otimes b^n\ps{n+1}h\ps{n+1}\otimes b^0b^1\ps{1}\cdots b^n\ps{1}h\ps{1}
\\
\qquad
=b^1\ps{2}\cdots b^n\ps{2}h\ps{2} \otimes\cdots\otimes \Delta\big(b^{j+1}\ps{j+1}\cdots b^n\ps{j+1}h\ps{j+1}\big) \otimes\cdots
\\
\qquad
\phantom{=}{}
\otimes b^n\ps{n+1}h\ps{n+1}\otimes b^0b^1\ps{1}\cdots b^n\ps{1}h\ps{1}
\\
\qquad
=s_j \circ \gamma_n\big(h\otimes b^0\otimes\cdots\otimes b^n\big).
\end{gather*}
Finally we observe that
\begin{gather*}
\gamma_n\circ t_n\big(h\otimes b^0\otimes\cdots\otimes b^n\big) =\gamma_n\big(b^n\ps{1}h\otimes b^n\ps{2}\otimes
b^0\otimes\cdots\otimes b^{n-1}\big)
\\
\qquad
=b^0\ps{2}\cdots b^{n-1}\ps{2}b^n\ps{1}\ps{2}h\ps{2} \otimes\cdots\otimes b^{n-1}\ps{n+1}b^n\ps{1}\ps{n}h\ps{n+1}\\
\qquad\hphantom{=}{}
\otimes
b^n\ps{2}b^0\ps{1}\cdots b^{n-1}\ps{1}b^n\ps{1}\ps{1}h\ps{1}
\\
\qquad
=b^0\ps{2}\cdots b^n\ps{2}h\ps{2} \otimes\cdots\otimes b^{n-1}\ps{n+1}b^n\ps{n+1}h\ps{n+1}\otimes
b^n\ps{n+2}b^0\ps{1}\cdots b^{n-1}\ps{1}b^n\ps{1}h\ps{1}
\\
\qquad
=b^0b^1\ps{1}\cdots b^n\ps{1}h\ps{2}\ps{2} \otimes\cdots\otimes b^{n-1}\ps{n}b^n\ps{n}h\ps{2}\ps{n+1}\otimes
b^n\ps{n+1}h\ps{3}S(h\ps{1})h\ps{2}\ps{1}
\\
\qquad
=b^0b^1\ps{1}\cdots b^n\ps{1}h\ps{1}\otimes b^1\ps{2}\cdots b^n\ps{2}h\ps{2} \otimes\cdots\otimes b^n\ps{n+1}h\ps{n+1}
\\
\qquad
=t_n \circ \gamma_n\big(h\otimes b^0\otimes\cdots\otimes b^n\big).  \tag*{\qed}
\end{gather*}
\renewcommand{\qed}{}
\end{proof}

\subsection{The classical case}

For any af\/f\/ine scheme~$X$ over a~f\/ield~$k$ we will denote by $\mathscr{O}(X)$ the~$k$-algebra of regular functions on~$X$.
Assume that~$G$ is a~linear algebraic group over a~f\/ield~$k$ and~$H$ its closed algebraic subgroup.
Assume that both are linearly reductive, hence all quotients below exist as af\/f\/ine varieties.

\subsubsection{The direct picture in the classical case}

The direct picture for ${\cal H}=\mathscr{O}(G)$, ${\cal C}=\mathscr{O}(H)$, ${\cal B}={\cal H}^{\co{\cal C}}=\mathscr{O}(G/H)$,
simplif\/ies a~lot because of the commutativity of ${\cal H}$.
Since ${\cal B}$ is central in ${\cal H}$, the cyclic tensor power $[{\cal H}^{\otimes_{{\cal B}} n+1}]_{{\cal B}}$
coincides with the usual tensor power ${\cal H}^{\otimes_{{\cal B}} n+1}$.
Moreover, the left SAYD ${\cal H}$-action on $\ad({\cal H})$ factors through the counit, and hence we can replace
$\ad({\cal H})$ by ${\cal H}$ regarded as a~trivial ${\cal H}$-module.
The right diagonal ${\cal H}$-action on the tensor power of ${\cal H}/{\cal I}$ comes from the multiplication in~$G$ of
elements of the subgroup~$H$.
Moreover, both sides become~$k$-algebras and the isomorphism~\eqref{aux-map-vp} becomes an algebra map
\begin{gather*}
\varphi_n\colon \  {\cal H}^{\otimes_{\cal B} n+1} \stackrel{\cong}{\longrightarrow} \big(({\cal H}/{\cal I})^{\otimes n+1}
\otimes_{\cal H} k\big)\otimes {\cal H},
\end{gather*}
which describes a~map of algebraic varieties
\begin{gather}
G^{\times_{G/H}\ n+1}\longleftarrow (H^{n+1}\times_{G}{\rm pt})\times G,
\nonumber
\\
(g, gh_{0}, gh_{0}h_{1},\ldots, gh_{0}\dots h_{n-1})\leftarrow\hspace{-0.18em}\shortmid (h_{0}, \ldots, h_{n}; g),
\label{str-pic-var-phi}
\end{gather}
with the inverse
\begin{gather}
G^{\times_{G/H}\ n+1}\longrightarrow (H^{n+1}\times_{G}{\rm pt})\times G,
\nonumber
\\
(g_{0}, \ldots, g_{n})\mapsto \big(g_{0}^{-1}g_{1}, g_{1}^{-1}g_{2}, \ldots, g_{n}^{-1}g_{0}; g_{0}\big).
\label{str-pic-var-psi}
\end{gather}
The left hand side is the f\/iber power of the canonical projection $G\rightarrow G/H$ consisting of tuples
$(g_{0},\ldots, g_{n})$ such that $g_{0}H=\cdots =g_{n}H$.
These varieties form a~cocyclic variety with the coface maps
\begin{gather*}
\delta_{i}\colon \  G^{\times_{G/H}\ n}\longrightarrow G^{\times_{G/H}\ n+1},  \qquad
\delta_{i}(g_{0}, \ldots, g_{n})=(g_{0}, \ldots, g_{i}, g_{i}, \ldots, g_{n}),
\end{gather*}
the codegeneracy maps
\begin{gather*}
\sigma_{i}\colon \ G^{\times_{G/H}\ n+1}\longrightarrow G^{\times_{G/H}\ n},  \qquad
\sigma_{i}(g_{0}, \ldots, g_{n})=(g_{0}, \ldots, g_{i-1},g_{i+1}, \ldots, g_{n}),
\end{gather*}
and the cocyclic map
\begin{gather*}
\tau_{n}\colon \  G^{\times_{G/H}\ n+1}\longrightarrow G^{\times_{G/H}\ n+1},  \qquad
\tau_{n}(g_{0}, g_{1},\ldots, g_{n})=(g_{1}, \ldots, g_{n}, g_{0}).
\end{gather*}
We note that this is the cyclic dual of the \v{C}ech nerve of the orbital map $G\rightarrow G/H$.

The right hand side is the cartesian product of the f\/iber product $H^{n+1}\times_{G}{\rm pt}$ of the map
$H^{n+1}\rightarrow G$, $(h_{0},\ldots, h_{n})\mapsto h_{0}\cdots h_{n}$ and the map ${\rm pt}\rightarrow G$,
$\star\mapsto e$, consisting of tuples $(h_{0},\ldots, h_{n})$ such that $h_{0}\cdots h_{n}=e$, and~$G$.
The cocyclic variety formed by these varieties has the coface maps
\begin{gather*}
\delta_{i}\colon \  \big(H^{n}\times_{G}{\rm pt}\big)\times G\longrightarrow \big(H^{n+1}\times_{G}{\rm pt}\big)\times G,
\\
\delta_{i}(h_{0}, \ldots, h_{n-1}; g)=(h_{0}, \ldots, h_{i-1}, e, h_{i}, \ldots, h_{n-1}; g),
\end{gather*}
the codegeneracy maps
\begin{gather*}
\sigma_{i}\colon \ \big(H^{n+1}\times_{G}{\rm pt}\big)\times G\longrightarrow \big(H^{n}\times_{G}{\rm pt}\big)\times G,
\\
\sigma_{i}(h_{0}, \ldots, h_{n}; g)=(h_{0}, \ldots, h_{i}h_{i+1}, \ldots, h_{n}; g),
\end{gather*}
and the cocyclic map
\begin{gather*}
\tau_{n}\colon \ \big(H^{n+1}\times_{G}{\rm pt}\big)\times G\longrightarrow \big(H^{n+1}\times_{G}{\rm pt}\big)\times G,
\\
\tau_{n}(h_{0}, h_{1},\ldots, h_{n}; g)=(h_{1}, \ldots, h_{n}, h_{0}; gh_{0}).
\end{gather*}
It is easy to check that the maps~\eqref{str-pic-var-phi}, \eqref{str-pic-var-psi} are mutually inverse isomorphisms of
cocyclic schemes.

\begin{Remark}
It is worth explaining the classical meaning of the isomorphism~\eqref{aux-Hoch-Tor}.
Let ${\cal H}$ be the Hopf algebra $\mathscr{O}(G)$ of regular functions on a~linear algebraic group~$G$ over
a~f\/ield~$k$ of characteristic zero, $\mathfrak{g}$ the Lie algebra of~$G$, and $\mathfrak{m}\subset\mathscr{O}(G)$ the
maximal ideal at the neutral element.

The left hand side of~\eqref{aux-Hoch-Tor} can be computed by the Hochschild--Kostant--Rosenberg
theorem~\cite{HochKostRose62} as the vector space of regular dif\/ferential forms on~$G$,
\begin{gather*}
{\rm HH}_\bullet({\cal H})=\Omega^\bullet(G).
\end{gather*}
As for the right hand side we use the commutativity of $\mathscr{O}(G)$.
Then the left $\mathscr{O}(G)$-module structure on $\ad(\mathscr{O}(G))$, \emph{a priori} coming from the
left-right SAYD-module structure, factors through the trivial action on~$k$ by the counit.
From geometric point of view, the counit is simply the evaluation at the neutral element of~$G$, corresponding to the
maximal ideal $\mathfrak{m}\subseteq \mathscr{O}(G)$.
Then we can use the identif\/ication $k={\mathscr O}(G)/\mathfrak{m}$, and the regularity of the group variety at the
neutral element, to use Serre's formula
\begin{gather*}
\Tor^{R}_{\bullet}(R/\mathfrak{m},
R/\mathfrak{m})=\Exterior^{\hspace{-0.1em}\bullet}_{R/\mathfrak{m}}\big(\mathfrak{m}/\mathfrak{m}^2\big),
\end{gather*}
where $R={\mathscr O}(G)_{\mathfrak{m}}$ is the local ring at the neutral element.

Using f\/inally the identif\/ication $\mathfrak{m}/\mathfrak{m}^2=\mathfrak{g}^{*}$, we obtain on the right hand side
of~\eqref{aux-Hoch-Tor},
\begin{gather*}
\Tor^{{\cal H}}_{\bullet}(k,{\cal H})=\Tor^{\mathscr{O}(G)}_{\bullet}(k,\mathscr{O}(G))
=\Tor^{\mathscr{O}(G)}_{\bullet}(k,k)\otimes\mathscr{O}(G)
=\Exterior^{\hspace{-0.1em}\bullet}\hspace{0.15em}\mathfrak{g}^{*}\otimes
\mathscr{O}(G).
\end{gather*}
Since both sides of~\eqref{aux-Hoch-Tor} are $\mathscr{O}(G)$-modules and the isomorphism is
$\mathscr{O}(G)$-linear,~\eqref{aux-Hoch-Tor} reads simply as triviality of the bundle of regular forms on a~linear
algebraic group~$G$
\begin{gather*}
\Omega^\bullet(G)=\Exterior^{\hspace{-0.1em}\bullet}\hspace{0.15em}\mathfrak{g}^{*}\otimes \mathscr{O}(G).
\end{gather*}
\end{Remark}

\subsubsection{The dual picture in the classical case}\label{section3.5.2}

The dual picture for ${\cal H}=\mathscr{O}(G)$, ${\cal C}=\mathscr{O}(H)$, ${\cal B}={\cal H}^{\co{\cal C}}=\mathscr{O}(G/H)$,
also can be made an isomorphism of cocyclic af\/f\/ine schemes over~$k$, however more interesting,
especially if ${\cal H}$ is not cocommutative, or equivalently, if~$G$ is not abelian.
Then the left-right SAYD module ${\rm coad}({\cal H})$ can be identif\/ied with $\mathscr{O}(\ad(G))$, where ${\rm
ad}(G)$ coincides with~$G$ as a~variety, the left ${\cal H}$-module structure comes from a~diagonal map of~$G$, and the
right ${\cal H}$-comodule structure comes from the right action of~$G$ on itself by conjugations.
Below it will be more convenient to pass to the left diagonal action of~$G$ on $\ad(G)\times (G/H)^{n+1}$ with the
use of the equivalent left~$G$-action on $\ad(G)$ by conjugations.
The isomorphism of the cyclic modules given by the maps
\begin{gather*}
{\rm coad}({\cal H})\Box_{{\cal H}}\big({\cal B}^{\otimes\ n+1}\big) \longrightarrow \big({\cal H}^{\Box_{{\cal C}}(n+1)}\big)^{{\cal C}}
\end{gather*}
reads now as an isomorphism of cocyclic varieties
\begin{gather*}
G^{n+1}/H^{n+1}\longrightarrow G\backslash \big(\ad(G)\times (G/H)^{n+1}\big).
\end{gather*}

On the left hand side we quotient by the right $H^{n+1}$-action
\begin{gather*}
G^{n+1}\times H^{n+1} \longrightarrow G^{n+1},
\\
(g_{0},\ldots, g_{n})\cdot (h_{0},\ldots, h_{n}) :=\big(h_{0}^{-1}g_{0}h_{1}, h_{1}^{-1}g_{1}h_{2},\ldots,
h_{n}^{-1}g_{n}h_{0}\big).
\end{gather*}
Note that this action is not free.
The stabilizer of $(g_{0},\ldots, g_{n})$ consists of the ($n+1$)-tuples $(h_{0},\ldots, h_{n})$ such that
\begin{gather}
h_{0}\in {\rm C}_{G}(g_{0}\cdots g_{n})\cap\bigcap_{i=0}^{n}(g_{0}\cdots g_{i})H (g_{0}\cdots g_{i})^{-1},
\nonumber\\
h_{i}=(g_{0}\cdots g_{i-1})h_{0} (g_{0}\cdots g_{i-1})^{-1}
\qquad
\text{for}
\quad
i=1,\ldots, n.\label{stab-LHS}
\end{gather}
The cocyclic structure on the left hand side is given by the coface operators
\begin{gather*}
\delta_{i}\colon \  G^{n+1}/H^{n+1} \longrightarrow G^{n+2}/H^{n+2},
\qquad
\delta_{i}[g_{0}, g_{1},\ldots, g_{n}] =[g_{0}, \ldots, g_{i-1}, e, g_{i}, \ldots, g_{n}],
\end{gather*}
the codegeneracy operators
\begin{gather*}
\sigma_{i}\colon \  G^{n+1}/H^{n+1} \longrightarrow G^{n}/H^{n},
\qquad
\sigma_{i}[g_{0}, g_{1},\ldots, g_{n}] =[g_{0}, \ldots, g_{i}g_{i+1}, \ldots, g_{n}],
\end{gather*}
and the cocyclic operator
\begin{gather*}
\tau_{n}\colon \  G^{n+1}/H^{n+1} \longrightarrow G^{n+1}/H^{n+1},
\qquad
\tau_{n}[g_{0}, g_{1},\ldots, g_{n}] =[g_{1}, \ldots, g_{n}, g_{0}].
\end{gather*}

On the right hand side we quotient by the left~$G$-action
\begin{gather*}
G\times \big(\ad(G)\times (G/H)^{n+1}\big) \longrightarrow \ad(G)\times (G/H)^{n+1},
\\
g\cdot(\widetilde{g}, g_{0}H,\ldots, g_{n}H) :=\big(g\widetilde{g}g^{-1}, g g_{0}H,\ldots, g g_{n}H\big).
\end{gather*}
Note that this action is not free as well.
The stabilizer of $(\widetilde{g}, g_{0}H,\ldots, g_{n}H)$ consists of the elements
\begin{gather}
\label{stab-RHS}
g   \in {\rm C}_{G}(\widetilde{g})\cap\bigcap_{i=0}^{n}g_{i}H g_{i}^{-1}.
\end{gather}
The cocyclic structure on the right hand side is given by the coface operators
\begin{gather}
\delta_{i}\colon \  G\backslash (\ad(G)\times (G/H)^{n+1}) \longrightarrow G\backslash \big(\ad(G)\times (G/H)^{n+2}\big),
\nonumber
\\
\delta_{i}[\widetilde{g}, g_{0}H,\ldots, g_{n}H] =\left[\widetilde{g}, g_{0}H,\ldots, g_{i}H, g_{i}H,\ldots,g_{n}H\right],
\label{cofaces-RHS}
\end{gather}
the codegeneracy operators
\begin{gather}
\sigma_{i}\colon \  G\backslash \big(\ad(G)\times (G/H)^{n+1}\big) \longrightarrow G\backslash \big(\ad(G)\times (G/H)^{n}\big),
\nonumber
\\
\sigma_{i}[\widetilde{g}, g_{0}H,\ldots, g_{n}H] =[\widetilde{g}, g_{0}H,\ldots,g_{i-1}H, g_{i+1}H, \ldots, g_{n}H],
\label{codegen-RHS}
\end{gather}
and the cocyclic operator
\begin{gather}
\tau_{n}\colon \  G\backslash \big(\ad(G)\times (G/H)^{n+1}\big) \longrightarrow G\backslash \big(\ad(G)\times (G/H)^{n+1}\big),
\nonumber
\\
\tau_{n}[\widetilde{g}, g_{0}H,\ldots, g_{n}H] =\left[\widetilde{g}, g_{1}H,\ldots, g_{n}H, \widetilde{g}g_{0}H\right].
\label{cyclic-op-RHS}
\end{gather}

It follows readily that the maps
\begin{gather}
G^{n+1}/H^{n+1}\longrightarrow G\backslash \big(\ad(G)\times (G/H)^{n+1}\big),
\nonumber
\\
[g_{0}, g_{1},\ldots, g_{n}]\mapsto \left[g_{0}g_{1}\cdots g_{n}; g_{0}H, g_{0}g_{1}H,\ldots, g_{0}g_{1}\cdots g_{n}H\right]
\label{cocyc-iso}
\end{gather}
and
\begin{gather}
G^{n+1}/H^{n+1}\longleftarrow G\backslash \big(\ad(G)\times (G/H)^{n+1}\big),
\nonumber
\\
\big[g_{n}^{-1}\widetilde{g}g_{0}, g_{0}^{-1}g_{1}, g_{1}^{-1}g_{2},\ldots, g_{n-1}^{-1}g_{n}\big]\leftarrow \shortmid
\left[\widetilde{g}; g_{0}H, g_{1}H,\ldots, g_{n}H\right]
\label{cocyc-iso-inv}
\end{gather}
are well def\/ined, mutually inverse, and they intertwine the coface, codegeneracy and cocyclic ope\-ra\-tors.
Quite unexpected and remarkable fact is that this isomorphism of af\/f\/ine varieties which identif\/ies orbits of dif\/ferent
groups acting on dif\/ferent varieties identif\/ies also their stabilizers, as is evident from the comparison
of~\eqref{stab-LHS} and~\eqref{stab-RHS}.
This suggests that a~more appropriate description should involve algebraic quotient stacks instead of varieties.

We note also that in~\eqref{cofaces-RHS},~\eqref{codegen-RHS} and~\eqref{cyclic-op-RHS}, the orbit $G/H$ can be replaced
with an arbitrary left~$G$-variety~$X$.
This produces a~cocyclic object given by the coface operators
\begin{gather}
\delta_{i}\colon \  G\backslash \big(\ad(G)\times X^{n+1}\big) \longrightarrow G\backslash \big(\ad(G)\times X^{n+2}\big),
\nonumber
\\
\delta_{i}[\widetilde{g}, x_{0},\ldots, x_{n}] =\left[\widetilde{g}, x_{0},\ldots, x_{i}, x_{i},\ldots, x_{n}\right],
\label{cofaces-X}
\end{gather}
the codegeneracy operators
\begin{gather}
\sigma_{i}\colon \ G\backslash \big(\ad(G)\times X^{n+1}\big) \longrightarrow G\backslash \big(\ad(G)\times X^{n}\big),
\nonumber
\\
\sigma_{i}[\widetilde{g}, x_{0},\ldots, x_{n}] =[\widetilde{g}, x_{0},\ldots,\widehat{x_{i}}, \ldots, x_{n}],
\label{codegen-X}
\end{gather}
and the cocyclic operator
\begin{gather}
\tau_{n}\colon \  G\backslash \big(\ad(G)\times X^{n+1}\big) \longrightarrow G\backslash \big(\ad(G)\times X^{n+1}\big),
\nonumber
\\
\tau_{n}[\widetilde{g}, x_{0},\ldots, x_{n}] = [\widetilde{g}, x_{1},\ldots, x_{n}, \widetilde{g} x_{0} ].
\label{cocyc-X}
\end{gather}
The periodic cyclic homology of its cyclic object of regular functions can be regarded as
a~$\mathbb{Z}/2$-graded~$G$-invariant cohomology of an af\/f\/ine~$G$-variety~$X$ with the equivariant system of
coef\/f\/i\-cients in a~vector space of regular functions on $\ad(G)$.

Such an \emph{ad-twisted cohomology} has been introduced in~\cite{Masz15}, and is denoted by ${\rm H}^{\bullet}_{\ad(G)}(X)$.
It is a~module over the algebra ${\rm H}^{\bullet}_{\ad(G)}({\rm pt})$ of ad-twisted cohomology of a~one
point~$G$-variety.
The latter is isomorphic to the algebra of regular class functions on~$G$ put in the even degree~\cite{Masz15}.
For a~trivial group~$G$ and a~smooth af\/f\/ine variety~$X$ its ad-twisted cohomology coincides with its
$\mathbb{Z}/2$-graded de Rham cohomology ${\rm H}^{\bullet}_{\ad(G)}(X)={\rm H}^{\bullet}_{\rm
dR}(X)$~\cite{Masz15}.

In~\cite{Masz15} an isomorphism similar to~\eqref{cocyc-iso}, \eqref{cocyc-iso-inv} (on the level of the Connes complex
computing the cyclic homology) was used in order to construct a~generalized character
\begin{gather*}
{\rm Rep}(G)\rightarrow {\rm H}^{0}_{\ad(G)}(G/H)
\end{gather*}
transforming into the Chern character of an associated vector bundle under a~Chern--Weil-like map, so called
\emph{strong Cartan connection}.

Moreover, for a~f\/inite~$G$ we mention brief\/ly the following relation between our construction and the \emph{Frobenius
reciprocity}.
Let us denote by $X=H\backslash \ad(H)$ and $Y=G\backslash \ad(G)$ the af\/f\/ine (f\/inite) varieties (over an
algebraically closed f\/ield~$k$ of characteristic not dividing the order of~$G$) of conjugacy classes, and by $f\colon
X\rightarrow Y$ the f\/inite \'{e}tale morphism induced by the containment $H<G$.
We then have two maps between the algebras $A={\mathscr O}(X)$ and $B={\mathscr O}(Y)$ of regular functions.
One is the algebra map $B\rightarrow A$ describing the morphism~$f$, and the other is a~$B$-linear map ${\rm Tr}_{f}\colon
A\rightarrow B$ def\/ined by means of the evaluation at points
\begin{gather*}
({\rm Tr}_{f}(a))(y)=\sum\limits_{f(x)=y}a(x).
\end{gather*}
Let $\chi\in A$ be a~character of a~representation of the subgroup~$H$, and denote by $\chi\uparrow^{G}_{H}\in B$ the
character of the induced representation.
Our point is that the Frobenius reciprocity can be rewritten as
\begin{gather}
\label{frob-rec}
\chi\uparrow^{G}_{H}={\rm Tr}_{f}(\chi),
\end{gather}
and the justif\/ication comes from the following canonical decomposition of~$f$ according to our construction
\begin{gather*}
H\backslash (\ad(H)\times H/H)\stackrel{i}{\rightarrow} G\backslash (\ad(G)\times
G/H)\stackrel{p}{\rightarrow} G\backslash (\ad(G)\times G/G).
\end{gather*}
Indeed, we have ${\rm Tr}_{f}={\rm Tr}_{p}\circ {\rm Tr}_{i}$,
\begin{gather*}
\xymatrix{A\ar[r]^-{{\rm Tr}_{i}} & B'\ar[r]^-{{\rm Tr}_{p}} &B},
\end{gather*}
where $B'={\mathscr O}(G\backslash (\ad(G)\times G/H))$, ${\rm Tr}_{i}$ is an extension by zero outside the image
of the closed-open immersion~$i$, and ${\rm Tr}_{p}$ is the usual trace map along the f\/ibers of the f\/inite \'{e}tale
covering~$p$.
Since explicitly
\begin{gather*}
{\rm Tr}_{i}(a)([\widetilde{g}, gH])=\begin{cases}
a\big(\big[g^{-1}\widetilde{g}g\big]\big), & {\rm if}\quad g^{-1}\widetilde{g}g\in H,
\\
0, &{\rm otherwise},
\end{cases}
\end{gather*}
and
\begin{gather*}
{\rm Tr}_{p}(b')([\widetilde{g}])=\sum\limits_{gH\in G/H}b'([\widetilde{g}, gH]),
\end{gather*}
we have
\begin{gather*}
{\rm Tr}_{f}(a)([\widetilde{g}])=\sum\limits_{gH\in G/H,\; g^{-1}\widetilde{g}g\in H} a\big(\big[g^{-1}\widetilde{g}g\big]\big),
\end{gather*}
which applied to~\eqref{frob-rec} gives the classical formulation of the Frobenius reciprocity
\begin{gather*}
\chi\uparrow^{G}_{H}([\widetilde{g}])=\sum\limits_{gH\in G/H,\; g^{-1}\widetilde{g}g\in H} \chi\big(\big[g^{-1}\widetilde{g}g\big]\big).
\end{gather*}
Note that in~\eqref{frob-rec} the left hand side depends only on the representation theory, while the right hand side
depends only on the geometry of the map~$f$ of conjugation classes.
This suggests that the Frobenius reciprocity in the form of~\eqref{frob-rec} should be read in view of the perfect
bilinear pairing provided by the evaluation of the algebra ${\mathscr O}(G\backslash \ad(G))$ of class functions
(with a~canonical basis consisting of the \emph{irreducible characters}) against the center ${\rm Z}(kG)$ of the group
algebra (with a~canonical basis consisting of \emph{class sums} corresponding to the conjugacy classes).
The latter can be understood as a~duality between the irreducible representations and the conjugacy classes.

Another interesting fact about the inverse isomorphism~\eqref{cocyc-iso-inv} is that it restricts to the collection of
varieties of \emph{extended quotients} in the sense of~\cite{AubeBaumPlymSoll14}.
The extended quotient $G\backslash\hspace{-0.23em}\backslash X$ of a~(say left)~$G$-variety~$X$ has been def\/ined
in~\cite{AubeBaumPlymSoll14} as a~usual quotient $G\backslash \widetilde{X}$ of a~subvariety $\widetilde{X}\subseteq
\ad(G)\times X$ of pairs $(\widetilde{g}, x)$ such that $\widetilde{g}x =x$ with respect to the left~$G$-action
$g(\widetilde{g}, x)=(g\widetilde{g}g^{-1}, gx)$.
The extended quotient replaces the orbit by the variety of the conjugacy classes of the stabilizer, and plays a~role in
the local Langlands program~\cite{AubeBaumPlymSoll14}.
For a~f\/inite~$G$ it was used to def\/ine the \emph{orbifold cohomology}~\cite{ChenRuan04} in terms of the \emph{inertia
orbifold}.
Usually atributed to Chern--Ruan who explored its new \emph{orbifold cup-product}, on the additive level it was in fact
introduced earlier in a~paper by Brylinski and Nistor, where the Chern--Ruan cohomology arises as the periodic cyclic
homology of the convolution algebra of the groupoid associated to the orbifold~\cite[Corollary~5.10(ii)]{BrylNist94}.
This is another evidence that cyclic homology of (quantum) stacks would be the most appropriate framework for
considering our duality.

It is easy to see that the operators~\eqref{cofaces-X}--\eqref{cocyc-X}
restrict to the collection
of extended quotients $G\backslash\hspace{-0.23em}\backslash (X^{n+1})\subseteq G\backslash (\ad(G)\times X^{n+1})$
with respect to the diagonal~$G$-action on the cartesian powers of a~left~$G$-variety~$X$ making it a~cocyclic variety.
Note that the~$G$-action on cartesian powers play a~role in the problem of \emph{inertia factors} in Clif\/ford
theory~\cite{PahlPles87}.

It is also easy to see that for the orbit $X=G/H$, the inverse isomorphism~\eqref{cocyc-iso-inv} restricts to
$G\backslash\hspace{-0.23em}\backslash ((G/H)^{n+1})\subseteq G\backslash (\ad(G)\times (G/H)^{n+1})$, and
transforms it into a~subvariety of $G^{n+1}/H^{n+1}$ consisting of the orbits of $(n+1)$-tuples $(g_{0}, g_{1},\ldots,
g_{n})$ such that $g_{i+1}\cdots g_{n}\cdot g_{0}\cdots g_{i}\in H$ for all $i=0,\ldots, n$.

\subsection*{Acknowledgements}
The authors would like to thank the anonymous referees for their constructive comments impro\-ving the paper.
The paper was partially supported by the NCN grant 2011/01/B/ST1/06474.
S.~S\"utl\"u would like to thank his former PhD advisor B.~Rangipour for drawing his attention to the homology of the
coalgebra-Galois extensions, Institut des Hautes \'Etudes Scientif\/iques (IHES) for the hospitality provided during part
of this work, and f\/inally the organizers of the conference ``From Poisson Brackets to Universal Quantum Symmetries'',
held at IMPAN, Warsaw, for the stimulating environment provided.

\pdfbookmark[1]{References}{ref}
\LastPageEnding


\begin{thebibliography}{99}
\footnotesize \itemsep=0pt

\bibitem{AubeBaumPlymSoll14}
Aubert A.M., Baum P., Plymen R., Solleveld M., Geometric structure in smooth
  dual and local {L}anglands conjecture, \href{http://dx.doi.org/10.1007/s11537-014-1267-x}{\textit{Jpn.~J. Math.}} \textbf{9}
  (2014), 99--136, \href{http://arxiv.org/abs/1211.0180}{arXiv:1211.0180}.

\bibitem{BergmanHausknecht-book}
Bergman G.M., Hausknecht A.O., Co-groups and co-rings in categories of
  associative rings, \href{http://dx.doi.org/10.1090/surv/045}{\textit{Mathematical Surveys and Monographs}}, Vol.~45,
  Amer. Math. Soc., Providence, RI, 1996.

\bibitem{Bich13}
Bichon J., Hochschild homology of {H}opf algebras and free {Y}etter--{D}rinfeld
  resolutions of the counit, \href{http://dx.doi.org/10.1112/S0010437X12000656}{\textit{Compos. Math.}} \textbf{149} (2013),
  658--678, \href{http://arxiv.org/abs/1204.0687}{arXiv:1204.0687}.

\bibitem{BoneCiccDabrTarl04}
Bonechi F., Ciccoli N., D\c{a}browski L., Tarlini M., Bijectivity of the
  canonical map for the non-commutative instanton bundle, \href{http://dx.doi.org/10.1016/j.geomphys.2003.09.007}{\textit{J.~Geom.
  Phys.}} \textbf{51} (2004), 71--81, \href{http://arxiv.org/abs/math.QA/0306114}{math.QA/0306114}.

\bibitem{BoneCiccTarl02}
Bonechi F., Ciccoli N., Tarlini M., Noncommutative instantons on the 4-sphere
  from quantum groups, \href{http://dx.doi.org/10.1007/s002200200618}{\textit{Comm. Math. Phys.}} \textbf{226} (2002),
  419--432, \href{http://arxiv.org/abs/math.QA/0012236}{math.QA/0012236}.

\bibitem{BrowZhan08}
Brown K.A., Zhang J.J., Dualising complexes and twisted {H}ochschild
  (co)homology for {N}oetherian {H}opf algebras, \href{http://dx.doi.org/10.1016/j.jalgebra.2007.03.050}{\textit{J.~Algebra}}
  \textbf{320} (2008), 1814--1850, \href{http://arxiv.org/abs/math.RA/0603732}{math.RA/0603732}.

\bibitem{BrylNist94}
Brylinski J.L., Nistor V., Cyclic cohomology of \'etale groupoids,
  \href{http://dx.doi.org/10.1007/BF00961407}{\textit{$K$-Theory}} \textbf{8} (1994), 341--365.

\bibitem{Brze96}
Brzezi{\'n}ski T., Quantum homogeneous spaces as quantum quotient spaces,
  \href{http://dx.doi.org/10.1063/1.531517}{\textit{J.~Math. Phys.}} \textbf{37} (1996), 2388--2399,
  \href{http://arxiv.org/abs/q-alg/9509015}{q-alg/9509015}.

\bibitem{Brze97}
Brzezi{\'n}ski T., Quantum homogeneous spaces and coalgebra bundles,
  \href{http://dx.doi.org/10.1016/S0034-4877(97)85914-9}{\textit{Rep. Math. Phys.}} \textbf{40} (1997), 179--185,
  \href{http://arxiv.org/abs/q-alg/9704015}{q-alg/9704015}.



\bibitem{BrzeHaja99}
Brzezi{\'n}ski T., Hajac P.M., Coalgebra extensions and algebra coextensions of
  {G}alois type, \href{http://dx.doi.org/10.1080/00927879908826498}{\textit{Comm. Algebra}} \textbf{27} (1999), 1347--1367,
  \href{http://arxiv.org/abs/q-alg/9708010}{q-alg/9708010}.

\bibitem{BrzeHaja09}
Brzezi{\'n}ski T., Hajac P.M., {G}alois-type extensions and equivariant
  projectivity, \href{http://arxiv.org/abs/0901.0141}{arXiv:0901.0141}.

\bibitem{CartEile56}
Cartan H., Eilenberg S., Homological algebra, Princeton University Press,
  Princeton, NJ, 1956.

\bibitem{ChenRuan04}
Chen W., Ruan Y., A new cohomology theory of orbifold, \href{http://dx.doi.org/10.1007/s00220-004-1089-4}{\textit{Comm. Math.
  Phys.}} \textbf{248} (2004), 1--31, \href{http://arxiv.org/abs/math.AG/0004129}{math.AG/0004129}.

\bibitem{CollHartThom09}
Collins B., H{\"a}rtel J., Thom A., Homology of free quantum groups,
  \href{http://dx.doi.org/10.1016/j.crma.2009.01.021}{\textit{C.~R.~Math. Acad. Sci. Paris}} \textbf{347} (2009), 271--276,
  \href{http://arxiv.org/abs/0903.1686}{arXiv:0903.1686}.

\bibitem{Conn83}
Connes A., Cohomologie cyclique et foncteurs {${\rm Ext}^n$},
  \textit{C.~R.~Acad. Sci. Paris S\'er.~I Math.} \textbf{296} (1983), 953--958.

\bibitem{ConnMosc98}
Connes A., Moscovici H., Hopf algebras, cyclic cohomology and the transverse
  index theorem, \href{http://dx.doi.org/10.1007/s002200050477}{\textit{Comm. Math. Phys.}} \textbf{198} (1998), 199--246,
  \href{http://arxiv.org/abs/math.DG/9806109}{math.DG/9806109}.

\bibitem{ConnMosc99}
Connes A., Moscovici H., Cyclic cohomology and {H}opf algebras, \href{http://dx.doi.org/10.1023/A:1007527510226}{\textit{Lett.
  Math. Phys.}} \textbf{48} (1999), 97--108, \href{http://arxiv.org/abs/math.QA/9904154}{math.QA/9904154}.

\bibitem{DijkKoor94}
Dijkhuizen M.S., Koornwinder T.H., Quantum homogeneous spaces, duality and
  quantum {$2$}-spheres, \href{http://dx.doi.org/10.1007/BF01278478}{\textit{Geom. Dedicata}} \textbf{52} (1994), 291--315.

\bibitem{Eckm86}
Eckmann B., Cyclic homology of groups and the {B}ass conjecture,
  \href{http://dx.doi.org/10.1007/BF02621911}{\textit{Comment. Math. Helv.}} \textbf{61} (1986), 193--202.

\bibitem{FengTsyg91}
Feng P., Tsygan B., Hochschild and cyclic homology of quantum groups,
  \href{http://dx.doi.org/10.1007/BF02099132}{\textit{Comm. Math. Phys.}} \textbf{140} (1991), 481--521.

\bibitem{Getz91-92}
Getzler E., Cartan homotopy formulas and the {G}auss--{M}anin connection in
  cyclic homology, in Quantum Deformations of Algebras and their
  Representations ({R}amat-{G}an, 1991/1992; {R}ehovot, 1991/1992),
  \textit{Israel Math. Conf. Proc.}, Vol.~7, Bar-Ilan University, Ramat Gan,
  1993, 65--78.

\bibitem{GinzKuma93}
Ginzburg V., Kumar S., Cohomology of quantum groups at roots of unity,
  \href{http://dx.doi.org/10.1215/S0012-7094-93-06909-8}{\textit{Duke Math.~J.}} \textbf{69} (1993), 179--198.

\bibitem{HadfKrah05}
Hadf\/ield T., Kr{\"a}hmer U., Twisted homology of quantum {${\rm SL}(2)$},
  \href{http://dx.doi.org/10.1007/s10977-005-3118-2}{\textit{$K$-Theory}} \textbf{34} (2005), 327--360, \href{http://arxiv.org/abs/math.QA/0405249}{math.QA/0405249}.

\bibitem{HajaKhalRangSomm04-II}
Hajac P.M., Khalkhali M., Rangipour B., Sommerh{\"a}user Y., Hopf-cyclic
  homology and cohomology with coef\/f\/icients, \href{http://dx.doi.org/10.1016/j.crma.2003.11.036}{\textit{C.~R.~Math. Acad. Sci.
  Paris}} \textbf{338} (2004), 667--672, \href{http://arxiv.org/abs/math.KT/0306288}{math.KT/0306288}.

\bibitem{HajaKhalRangSomm04-I}
Hajac P.M., Khalkhali M., Rangipour B., Sommerh{\"a}user Y., Stable
  anti-{Y}etter--{D}rinfeld modules, \href{http://dx.doi.org/10.1016/j.crma.2003.11.037}{\textit{C.~R.~Math. Acad. Sci. Paris}}
  \textbf{338} (2004), 587--590, \href{http://arxiv.org/abs/math.QA/0405005}{math.QA/0405005}.

\bibitem{HochKostRose62}
Hochschild G., Kostant B., Rosenberg A., Dif\/ferential forms on regular af\/f\/ine
  algebras, \href{http://dx.doi.org/10.1090/S0002-9947-1962-0142598-8}{\textit{Trans. Amer. Math. Soc.}} \textbf{102} (1962), 383--408.

\bibitem{JaraStef06}
Jara P., {\c{S}}tefan D., Hopf-cyclic homology and relative cyclic homology of
  {H}opf--{G}alois extensions, \href{http://dx.doi.org/10.1017/S0024611506015772}{\textit{Proc. London Math. Soc.}} \textbf{93}
  (2006), 138--174.

\bibitem{Kadi89}
Kadison L., A relative cyclic cohomology theory useful for computations,
  \textit{C.~R.~Acad. Sci. Paris S\'er.~I Math.} \textbf{308} (1989), 569--573.

\bibitem{Kadi92}
Kadison L., Cyclic homology of triangular matrix algebras, in Topology {H}awaii
  ({H}onolulu, {HI}, 1990), World Sci. Publ., River Edge, NJ, 1992, 137--148.

\bibitem{Kadi14}
Kadison L., Hopf subalgebras and tensor powers of generalized permutation
  modules, \href{http://dx.doi.org/10.1016/j.jpaa.2013.06.008}{\textit{J.~Pure Appl. Algebra}} \textbf{218} (2014), 367--380,
  \href{http://arxiv.org/abs/1210.3178}{arXiv:1210.3178}.

\bibitem{KhalRang05}
Khalkhali M., Rangipour B., A note on cyclic duality and {H}opf algebras,
  \href{http://dx.doi.org/10.1081/AGB-200051130}{\textit{Comm. Algebra}} \textbf{33} (2005), 763--773,
  \href{http://arxiv.org/abs/math.KT/0310088}{math.KT/0310088}.

\bibitem{Lu93}
Lu J.H., Moment maps at the quantum level, \href{http://dx.doi.org/10.1007/BF02099767}{\textit{Comm. Math. Phys.}}
  \textbf{157} (1993), 389--404.

\bibitem{Masz15}
Maszczyk T., {F}eynman integral, {C}hern character and duality, {i}n
  preparation.

\bibitem{PahlPles87}
Pahlings H., Plesken W., Group actions on {C}artesian powers with applications
  to representation theory, \href{http://dx.doi.org/10.1515/crll.1987.380.178}{\textit{J.~Reine Angew. Math.}} \textbf{380} (1987),
  178--195.

\bibitem{Podl87}
Podle{\'s} P., Quantum spheres, \href{http://dx.doi.org/10.1007/BF00416848}{\textit{Lett. Math. Phys.}} \textbf{14} (1987),
  193--202.

\bibitem{Scha92}
Schafer J.A., Relative cyclic homology and the {B}ass conjecture,
  \href{http://dx.doi.org/10.1007/BF02566497}{\textit{Comment. Math. Helv.}} \textbf{67} (1992), 214--225.

\bibitem{Scha98-II}
Schauenburg P., Galois correspondences for {H}opf bi-{G}alois extensions,
  \href{http://dx.doi.org/10.1006/jabr.1997.7237}{\textit{J.~Algebra}} \textbf{201} (1998), 53--70.

\bibitem{Schn93}
Schneider H.J., Some remarks on exact sequences of quantum groups,
  \href{http://dx.doi.org/10.1080/00927879308824733}{\textit{Comm. Algebra}} \textbf{21} (1993), 3337--3357.

\bibitem{Take72}
Takeuchi M., A correspondence between {H}opf ideals and sub-{H}opf algebras,
  \href{http://dx.doi.org/10.1007/BF01579722}{\textit{Manuscripta Math.}} \textbf{7} (1972), 251--270.

\bibitem{VanOZhan94}
Van~Oystaeyen F., Zhang Y., Galois-type correspondences for {H}opf--{G}alois
  extensions, \href{http://dx.doi.org/10.1007/BF00960864}{\textit{$K$-Theory}} \textbf{8} (1994), 257--269.

\end{thebibliography}
\end{document}